\documentclass[11pt]{article}
\usepackage{graphicx}
\usepackage{amsmath,amssymb,amsthm,amsfonts}
\usepackage{amssymb}
\newtheorem{thm}{Theorem}[section]

\newtheorem{lem}[thm]{Lemma}
\newtheorem{prop}[thm]{Proposition}

\newtheorem{rem}{Remark}[section]
\theoremstyle{definition}
\newtheorem{defn}{Definition}[section]
\usepackage{appendix}
\newtheorem{ex}{Example}[section]

\numberwithin{equation}{section}

\DeclareMathSymbol{\C}{\mathalpha}{AMSb}{"43}

\newcommand{\eps}{\varepsilon}

\newcommand{\lam}{\lambda}
\newcommand{\alp}{\alpha}

\newcommand{\R}{{\mathbb{R}}}
\newcommand{\h}{{\mathcal{H}}}
\newcommand{\inte}{\int_{\mathbb{R}^2}}
\newcommand{\intr}{\int_0^\infty}
\newcommand{\intB}{\int _{B_\delta (x_{2,k})}}
\newcommand{\intPB}{\int _{\partial B_\delta (x_{2,k})}}

\newcommand{\bsub}{\begin{subequations}}
\newcommand{\esub}{\end{subequations}$\!$}

\begin{document}

\title{Local Uniqueness and Refined Spike Profiles of Ground States for Two-Dimensional Attractive Bose-Einstein Condensates}

\author{Yujin Guo\thanks{Wuhan Institute of Physics and Mathematics,
    Chinese Academy of Sciences, P.O. Box 71010, Wuhan 430071,
    P. R. China.  Email: \texttt{yjguo@wipm.ac.cn}. Y. J.  Guo is partially supported by NSFC grants 11322104 and 11671394.
    },
    \, Changshou Lin\thanks{Taida Institute of Mathematical Sciences, National Taiwan University,  Taipei 10617, Taiwan.  Email: \texttt{cslin@math.ntu.edu.tw}. },
\, and\, Juncheng Wei\thanks{Department of Mathematics, University of British Columbia, Vancouver, BC V6T 1Z2,
    Canada.  Email: \texttt{jcwei@math.ubc.ca}. J.C. Wei is partially supported by NSERC of Canada.
    }
    }

\date{\today}

\smallbreak \maketitle

\begin{abstract}
We consider ground states of two-dimensional Bose-Einstein condensates in a trap with attractive interactions, which can be described equivalently by positive minimizers of the $L^2-$critical constraint Gross-Pitaevskii energy functional. It is known that ground states exist if and only if $a< a^*:= \|w\|_2^2$, where $a$ denotes the  interaction strength and $w$ is the unique positive solution of $\Delta w-w+w^3=0$ in $\R^2$. In this paper, we prove the local uniqueness and refined spike profiles of ground states as $a\nearrow a^*$, provided that  the trapping potential $h(x)$ is  homogeneous and $H(y)=\inte h(x+y)w^2(x)dx$ admits a unique and non-degenerate critical point.
\end{abstract}

\vskip 0.1truein


\noindent {\it Keywords:} Bose-Einstein condensation; spike profiles; local uniqueness; Pohozaev identity.

\vskip 0.2truein

\section{Introduction}

The phenomenon of Bose-Einstein condensation (BEC) has been investigated intensively since its first realization in cold atomic gases, see   \cite{Anderson,D} and references therein. In these experiments, a large number
of (bosonic) atoms are confined to a trap and cooled to very low temperatures. Condensation of a large fraction
of particles into the same one-particle state is observed below a  critical temperature.
These Bose-Einstein condensates display various interesting quantum phenomena,
such as the critical-mass collapse, the superfluidity and the appearance of quantized vortices in rotating
traps (e.g.\cite{D}). Specially, if the force between the atoms in the condensates is attractive, the system collapses as soon as the particle number increases
beyond a critical value, see, e.g., \cite{Hulet3} or
\cite[Sec.~III.B]{D}.

Bose-Einstein condensates (BEC)  of a dilute gas with attractive interactions in $\R^2$ can be described (\cite{Bao2,D,GS}) by the following Gross-Pitaevskii (GP) energy functional
\begin{equation}\label{f}
  E_a(u):=\int_{\R ^2}\Big(|\nabla
  u|^2+V(x)|u|^2\Big)dx-\frac{a}{2}\inte|u|^4dx,
\end{equation}
where $a>0$ describes the strength of the attractive interactions, and $V(x)\ge 0$ denotes the trapping potential satisfying $\lim_{|x|\to\infty} V(x) = \infty$. As addressed recently in \cite{GS,GWZZ}, ground states of attractive BEC in $\R^2$ can be described by the constraint minimizers of the GP energy
\begin{equation}\label{eq1.3}
e(a):=\inf _{\{u\in \h, \, \|u\|^2_2=1 \} } E_a(u),
\end{equation}
where the space $\h$ is defined by
\begin{equation}
   \h := \Big \{u\in  H^1(\R ^2):\ \int _{\R ^2}  V(x)|u(x)|^2 dx<\infty \Big\}. \label{H}
\end{equation}
The minimization problem $e(a)$ was analyzed recently in \cite{Bao2,GS,GWZZ,GZZ,Z} and references therein. Existing results show that $e(a)$ is an $L^2-$critical constraint variational problem. Actually, it was shown in \cite{Bao2,GS} that $e(a)$ admits minimizers if and only if $a<a^*:=\|w\|^2_2$,  where $w=w(|x|)$ is the unique (up to translations) radial positive solution (cf. \cite{GNN,mcleod,K})  of the following nonlinear scalar
field equation
\begin{equation}
\Delta w-w+w^3=0\  \mbox{  in } \  \R^2,\  \mbox{ where }\  w\in
H^1(\R ^2)\,.  \label{Kwong}
\end{equation}
It turns out that the existence and nonexistence of minimizers for $e(a)$ are well connected with the following Gagliardo-Nirenberg inequality
\begin{equation}\label{GNineq}
\inte |u(x)|^4 dx\le \frac 2 {\|w\|_2^{2}} \inte |\nabla u(x) |^2dx
\inte |u(x)|^2dx ,\   \ \forall u \in H^1(\R ^2)\,,
\end{equation}
where the equality is attained at $w$ (cf. \cite{W}).

Since $E_a(u)\geq E_a(|u|)$ for any $u\in \mathcal{H}$, any minimizer $u_a$ of $e(a)$ must be either non-negative or non-positive, and  it   satisfies
the Euler-Lagrange equation
\begin{equation}
-\Delta u_a  +V(x)u_a  =\mu_a u_a +au ^3 _a \quad \  \text{in\,  $\R^2$},
\label{A:EL}
\end{equation}
where $\mu _a\in \R$ is a suitable Lagrange multiplier.  Thus, by applying the maximum principle to the equation (\ref{A:EL}), any minimizer $u_a$ of $e(a)$ is further either negative or positive. Therefore, without loss of generality one can restrict the
minimizations of $e(a)$ to positive functions.  {\em In this paper positive minimizers of $e(a)$ are called ground states of attractive BEC.}
Applying energy estimates and blow-up analysis, the spike profiles of positive minimizers for $e(a)$ as $a\nearrow a^*$ were recently discussed in \cite{GS,GWZZ,GZZ} under different types of potentials $V(x)$, see our Proposition \ref{2:prop} for some related results. In spite of these facts, it remains open to discuss the refined spike profiles of positive minimizers. On the other hand, the local uniqueness of positive minimizers for $e(a)$ as $a. e.$ $a\nearrow a^*$ was also proved \cite{GWZZ} by the ODE argument, for the case where $V(r)=V(|x|)$ is radially symmetric and satisfies $V'(r)\ge 0$, see Corollary 1.1 in \cite{GWZZ} for details. Here the locality of uniqueness means that   $a$ is near $ a^*$. It is therefore natural to ask whether such local uniqueness still holds for the case where $V(x)$ is not radially symmetric. We should remark that all these results mentioned above were obtained mainly by analyzing the variational structures of the minimization problem $e(a)$, instead of discussing the PDE properties of the associated elliptic equation (\ref{A:EL}).

By investigating thoroughly the associated equation (\ref{A:EL}), the main purpose of this paper is to derive the refined spike profiles of positive minimizers for $e(a)$ as $a\nearrow a^*$, and  extend the above local uniqueness to the cases of non-symmetric  potentials $V(x)$ as well. Throughout the whole paper, we shall consider the trapping potential $V(x)$ satisfying $\lim_{|x|\to\infty} V(x) = \infty$ in the class of homogeneous functions, for which we define
\begin{defn}
$h(x)\ge 0$ in $\R^2$ is homogeneous of degree $p\in\R^+$ (about the origin), if there exists some $p>0$ such that
 \begin{equation}\label{1:V}
h(tx)=t^ph(x)\ \, \mbox{in}\ \, \R^2 \ \mbox{for any}\ t>0.
\end{equation}
\end{defn}
\noindent Following \cite[Remark 3.2]{Grossi}, the above definition implies that the homogeneous function $h(x)\in C(\R^2)$ of degree $p>0$ satisfies
 \begin{equation}\label{1:Vh}
0\le h(x)\le C|x|^p\,\ \mbox{in}\,\ \R^2,
\end{equation}
where $C>0$ denotes the maximum of $h(x)$ on $\partial B_1(0)$. Moreover, since we assume that $\lim_{|x|\to\infty} h(x) = \infty$, $x=0$ is the unique minimum point of $h(x)$. Additionally, we often need to assume that $V(x)=h(x)\in C^2(\R^2)$ satisfies
\begin{equation}\label{1:HH}
y_0 \,\ \text{is the unique critical point of}\,\ H(y)=\inte h(x+y)w^2(x)dx.
\end{equation}
The following example shows that for some non-symmetric potentials $h(x)$, $H(y)$ admits a unique critical point  $y_0$, where $y_0$ satisfies $y_0\not =0$ and is {\em non-degenerate} in the sense that
\begin{equation}\label{1:HK}
 \det \Big(\frac{\partial ^2H(y_0)}{\partial x_i\partial x_j}\Big)\not =0,\,\ \text{where}\,\ i,\, j=1,\, 2.
\end{equation}

\begin{ex}\label{1:example}  Suppose that the potential $h(x)$ satisfies
\begin{equation}\label{1:PPP}
 h(x)=|x|^p\big[1+\delta h_0(\theta)\big]\ge 0, \,\ \text{where}\,\ p\ge 2 \,\ \text{and}\,\ \delta \in \R,
\end{equation}
where $h_0(\theta)\in C^2([0,2\pi ])$ satisfies
\begin{equation}\label{1:PPP-1}
 \Big(\int _{0}^{2\pi}h_0(\theta)\cos \theta d\theta\Big)^2+\Big(\int _{0}^{2\pi}h_0(\theta)\sin \theta d\theta\Big)^2>0.
\end{equation}
One can check from (\ref{1:PPP-1}) that if $|\delta|\ge 0$ is small enough, then $H(y)$ admits a unique critical point $y_0=-\delta \hat y_0\in \R^2$, where $\hat y_0$ satisfies
\begin{equation}\label{1:PPP-2}
\hat y_0\, \sim\, \Big(C_1\int _{0}^{2\pi}h_0(\theta)\cos \theta d\theta,\, C_2\int _{0}^{2\pi}h_0(\theta)\sin \theta d\theta\Big)\not =(0,0)\,\ \text{as}\,\  \delta \to 0
\end{equation}
for some positive constants $C_1$ and $C_2$ depending only on $w$ and $p$. Furthermore, if $|\delta|\ge 0$ is small enough, then $\det \Big(\frac{\partial ^2H(y_0)}{\partial x_i\partial x_j}\Big)>0$, which implies that the unique critical point $y_0$ of $H(y)$ is non-degenerate.
\end{ex}

Our first main result is concerned with the following local uniqueness as $a\nearrow a^*$, which holds for some non-symmetric homogeneous potentials $h(x)$ in view of Example 1.1.

\begin{thm}\label{1:thmA}
Suppose $V(x)=h(x)\in C^2(\R^2)$ is homogeneous of degree $p\ge 2$, where  $\lim_{|x|\to\infty} h(x) = \infty$, and satisfies
\begin{equation}\label{1:H}
y_0 \,\ \text{is the unique and non-degenerate critical point of}\,\ H(y)=\inte h(x+y)w^2(x)dx.
\end{equation}
Then there exists a unique positive minimizer for $e(a)$ as $a\nearrow a^*$.
\end{thm}

%


The local uniqueness of Theorem \ref{1:thmA} means that  positive minimizers of $e(a)$ must be unique as $a$ is near $ a^*$.
It is possible to extend Theorem \ref{1:thmA} to  more general potentials $V(x)=g(x)h(x)$ for a class of functions $g(x)$, which is however beyond the discussion ranges of the present paper.
We also remark that the proof of Theorem \ref{1:thmA} is more involved for the case where $y_0\not =0$ occurs in (\ref{1:H}).
Our proof of such local uniqueness is motivated by \cite{Cao,Deng,Grossi}. Roughly speaking, as derived in Proposition \ref{2:prop} we shall first obtain some fundamental estimates on the spike behavior of positive minimizers. Under the non-degeneracy assumption of (\ref{1:H}), the local uniqueness is then proved in Subsection 2.1 by establishing various types of local Pohozaev identities.

The proof of Theorem \ref{1:thmA} shows that if one considers the local uniqueness of Theorem \ref{1:thmA} in other dimensional cases, where $\R^2$ is replaced by $\R^d$  and $u^4$ is replaced by $u^{2+\frac{4}{d}}$ for $d\not =2$, the fundamental estimates of Proposition \ref{2:prop} are not enough. Therefore, in the following we address the refined  spike behavior of positive minimizers under the assumption (\ref{1:H}). To introduce our second main result, for convenience we next denote
\begin{equation}\label{1a:li}
  \lambda  _0= \left(\frac{p}{2}\int_{\R^2} h(x+y_0) w^2(x)dx   \right)^{\frac 1{2+p}},
\end{equation}
where $y_0\in \R^2$ is given by (\ref{1:H}), and
\[
\psi (x)=\varphi (x)-\frac{C^*}{2}\big[w(x)+x\cdot \nabla w(x)\big],
\]
where $\varphi (x)\in C^2(\R^2)\cap L^\infty(\R^2)$  is the unique solution of
     \begin{equation}\label{thm1:eqn1}
 \nabla \varphi (0)=0\ \, \mbox{and}\,\ \big[-\Delta +(1-3w^2)\big]\varphi (x)=-\displaystyle \frac{2w^3}{\inte w^4}-\displaystyle\frac{2h(x+y_0)w}{p\inte h(x+y_0)w^2}
 \,\ \mbox{in}\,\ \R^2,
   \end{equation}
and the nonzero constant $C^*$ is given by
\[
C^*= \frac{2}{2+p}\Big(2\displaystyle\inte w\psi_{3}+\displaystyle\inte \varphi^2\Big)
\]
with $\psi_{3}\in C^2(\R^2)\cap L^\infty(\R^2)$ being the unique solution of (\ref{lem3.2:2}). Using above notations, we shall derive the following theorem.

\begin{thm}\label{1:thmB}
Suppose $V(x)=h(x)\in C^2(\R^2)$ is homogeneous of degree $p\ge 2$, where  $\lim_{|x|\to\infty} h(x) = \infty$, and  satisfies
(\ref{1:H}) for some $y_0\in \R^2$. If $u_a$ is a positive minimizer of $e(a)$ as $a\nearrow a^*$, then we have
\begin{equation}\arraycolsep=1.5pt\begin{array}{lll}
u_a(x)&=& \displaystyle\frac { \lambda  _0}{\|w\|_2} \Big\{  \displaystyle\frac{1}{(a^*-a)^\frac{1}{2+p}} w\Big(\frac{\lam _0 (x-x_a)}{(a^*-a)^\frac{1}{2+p}}\Big)+\displaystyle(a^*-a)^\frac{1+p}{2+p}\psi\Big(\frac{\lam _0 (x-x_a)}{(a^*-a)^\frac{1}{2+p}}\Big)\\[4mm]
&&+\displaystyle(a^*-a)^\frac{3+2p}{2+p} \phi_0\Big(\frac{\lam  _0(x-x_a)}{(a^*-a)^\frac{1}{2+p}}\Big) \Big\}  +o\big((a^*-a)^\frac{3+2p}{2+p}\big)   \ \ \mbox{as} \ \  a \nearrow a^*
\end{array}\label{1a:cona}
\end{equation}
uniformly in $\R^2$ for some function $\phi_0\in C^2(\R^2)\cap L^\infty(\R^2)$, where $ x_a $ is the unique maximum point of $u_a$ satisfying
\begin{equation}
 \Big|\frac{\lam x_a }{(a^*-a)^\frac{1}{2+p}}-y_0\Big|=(a^*-a)O(|y^0|) \ \ as \  \ a  \nearrow a^*
\label{1a:conb}
\end{equation}
for some $y^0\in\R^2$.
\end{thm}

%



Theorem \ref{1:thmB} is derived directly from  Theorem \ref{1:thmA} and Theorem    \ref{thm3.5} in Section 3  with more details, where $\phi_0\in C^2(\R^2)\cap L^\infty(\R^2)$ is given explicitly. In Section 4 we shall extend the refined  spike behavior of Theorem \ref{1:thmB} to more general potentials $V(x)=g(x)h(x)$, where
$h(-x)=h(x)$ is homogeneous and satisfies (\ref{1:H}) and $0\le C\le g(x)\le \frac{1}{C}$ holds in $\R^2$, see Theorem \ref{thm4.5} for details. To establish Theorem \ref{1:thmB} and Theorem \ref{thm4.5}, our Proposition \ref{2:prop} shows that the arguments of \cite{GS,GWZZ,GZZ} give the leading expansion terms of the minimizer $u_a$ and the associated  Lagrange multiplier $\mu _a$ satisfying (\ref{A:EL}) as well. In order to get (\ref{1a:cona}) for the rest terms of $u_a$, the difficulty is to obtain the more precise estimate of $\mu _a$, which is overcome by the very delicate analysis of the associated equation (\ref{A:EL}), together with the constraint condition of $u_a$.

This paper is organized as follows: In Section 2 we shall prove Theorem \ref{1:thmA}  on the local uniqueness of positive minimizers. Section 3 is concerned with proving Theorem \ref{1:thmB} on the refined spike profiles of positive minimizers for $e(a)$ as $a\nearrow a^*$. The main aim of Section 4 is to derive Theorem \ref{thm4.5}, which extends the refined  spike behavior of Theorem \ref{1:thmB} to more general potentials $V(x)=g(x)h(x)$. We shall leave the proof of Lemma \ref{lem3.3} to Appendix A.

\section{Local Uniqueness of Positive Minimizers}

This section is devoted to the proof of  Theorem \ref{1:thmA} on the local uniqueness of positive minimizers.
Towards this purpose, we need some estimates of positive minimizers for $e(a)$ as $a\nearrow a^*$, which hold essentially for more general potential $V(x)\in C^2 (\R^2)$ satisfying
\begin{equation}\arraycolsep=1.5pt\begin{array}{lll}
&& V(x)=g(x)h(x),\, \mbox{where $0<C\le g(x)\le \frac{1}{C}$ in $\R^2$ and $h(x)$ is homogeneous of}\\[2mm]
   &&\mbox{degree $p\ge 2$.}
 \end{array}
\label{2:V1}
\end{equation}
For convenience, we always denote $\{u_k\}$ to be a
positive minimizer sequence of $e(a_k)$ with $a_k\nearrow a^*$ as $k\to\infty$, and define
\begin{equation}\label{def:li}
  \lambda = \left( \frac{pg(0)}{2}  \int_{\R^2} h(x+y_0) w^2(x)dx   \right)^{\frac 1{2+p}},
\end{equation}
where $V(x)=g(x)h(x)$ is assumed to satisfy (\ref{2:V1}) with $p\ge 2$ and  $y_0\in \R^2$ is given by (\ref{1:HH}). Recall  from (\ref{Kwong})  that $w(|x|)$  satisfies
\begin{equation}\label{1:id}
\inte |\nabla w |^2dx  =\inte |w| ^2dx=\frac{1}{2}\inte |w| ^4dx ,
\end{equation} see also Lemma 8.1.2 in \cite{C}.  Moreover, it follows from \cite[Prop.~4.1]{GNN} that $w$ admits the following exponential decay
 \begin{equation}
w(x) \, , \ |\nabla w(x)| = O(|x|^{-\frac{1}{2}}e^{-|x|}) \quad
\text{as \ $|x|\to \infty$.}  \label{1:exp}
\end{equation}

\begin{prop}\label{2:prop}
Suppose $V(x)=g(x)h(x)\in C^2 (\R^2)$ satisfies $\lim_{|x|\to\infty} V(x) = \infty$ and (\ref{2:V1}), and assume (\ref{1:HH}) holds for some $y_0\in \R^2$. Then there exist a subsequence, still denoted by $\{a_k\}$, of $\{a_k\}$ and $\{x_k\}\subset\R^2$ such that
\begin{enumerate}
\item [(I).] The subsequence $\{u_k\}$ satisfies
\begin{equation}
(a^*-a_k)^{\frac{1}{2+p}} u_k\left(x_k + x  (a^*-a_k)^{\frac{1}{2+p}}\right) \to  \frac { \lambda w(\lambda x)}{\|w\|_2}    \,\ \mbox{as} \,\ k\to\infty \label{2:cona}
\end{equation}
uniformly in $\R^2$, and $ x_k $ is the unique maximum point of $u_k$ satisfying
\begin{equation}
\lim_{k\to \infty} \frac{\lam x_k}{(a^*-a_k)^{\frac{1}{2+p}}}=y_0,
\label{2:conb}
\end{equation}
where $y_0\in \R^2$ is the same as that of (\ref{1:HH}).
Moreover, $u_k$ satisfies
 \begin{equation}
(a^*-a_k)^{\frac{1}{2+p}} u_k\left(x_k + x  (a^*-a_k)^{\frac{1}{2+p}}\right)\le Ce^{-\frac{\lam}{2}|x|} \,\ \text{in\,\ $\R^2$},
\label{2:conexp}
\end{equation}
where the constant $C>0$ is independent of $k$.
\item [(II).] The energy $e(a_k)$ satisfies
\begin{equation}\label{2:conc}
\lim_{k\to \infty} \frac{e(a_k)}{(a^*-a_k)^{p/(2+p)}} = \frac{\lambda^2}{a^*} \frac{p+2}{p}.
\end{equation}
\end{enumerate}
\end{prop}

\noindent\textbf{Proof.}
Since the proof of Proposition \ref{2:prop} is similar to those in \cite{GS,GWZZ,GZZ}, which handle (\ref{f}) with different potentials $V(x)$, we shall briefly sketch the structure of the proof.

If $V(x)\in C^2(\R^2)$ satisfies (\ref{2:V1}) with $p\ge 2$, we note that $h(x)\ge 0$ satisfies (\ref{1:Vh}). Take the test function
\[
u_\tau (x) =  A_{ \tau} \frac{\tau}{\|w\|_2} \varphi( x  ) w(\tau x),
\]
where the nonnegative cut-off function $\varphi \in C_0^\infty(\R^2)$ satisfies $0\le \varphi(x) \le 1$ in $\R^2$, and $A_{ \tau}>0$ is chosen so that $\inte u_\tau (x)^2dx=1$. The same proof of Lemma 3 in \cite{GS} then yields that
\begin{equation}
 e(a)\le C(a^*-a)^{\frac{p}{p+2}}\ \ \mbox{for}\ \ 0\leq a<  a^*,
\label{2:con:1}
\end{equation}
where the constant $C>0$ is independent of $a$. By (\ref{2:con:1}), we can follow Lemma 4 in \cite{GS} to derive that
there exists a positive constant  $K$, independent of $a$, such that
\begin{equation}
\inte |u_a(x)|^4dx \le \frac{1}{K}(a^*-a)^{-\frac{2}{p+2}}\ \ \mbox{for}\ \ 0\leq a<  a^*,
\label{2:con:2}
\end{equation}
where $u_a>0$ is any minimizer of $e(a)$. Applying (\ref{2:con:1}) and (\ref{2:con:2}), a proof similar to that of Theorem 2.1 in \cite{GZZ} then gives that there exist two positive constants $m<M$, independent of $a$, such that
\begin{equation}
m(a^*-a)^{\frac{p}{p+2}}\le e(a)\le M(a^*-a)^{\frac{p}{p+2}}\ \ \mbox{for}\ \ 0\leq a< a^*.
\label{2:con:3}
\end{equation}
Based on (\ref{2:con:3}), similar to Theorems 1.2 and 1.3 in \cite{GZZ}, one can further deduce that there exist a subsequence  (still denoted by $\{a_k\}$) of $\{a_k\}$ and $\{x_k\}\subset\R^2$, where $a_k\nearrow a^*$ as $k\to\infty$, such that (\ref{2:conexp}) and (\ref{2:conc}) hold,  and
\begin{equation}
(a^*-a_k)^{\frac{1}{2+p}} u_k\left(x_k + x  (a^*-a_k)^{\frac{1}{2+p}}\right) \to  \frac { \lambda w(\lambda x)}{\|w\|_2}  \,  \,\ \mbox{strongly \ in} \,\ H^1(\R^2) \label{2:con:cona}
\end{equation}
as $k\to\infty$, where $ x_k $ is the unique maximum point of $u_k$. Finally, since $w$ decays exponentially, the standard elliptic regularity theory applied to (\ref{2:con:cona}) yields that (\ref{2:cona}) holds uniformly in $\R^2$ (e.g. Lemma 4.9 in \cite{M} for similar arguments).

We finally follow (\ref{1:HH}) and  (\ref{2:cona}) to derive the estimate  (\ref{2:conb}). Following (\ref{2:cona}), we define
\[
\bar u_k(x):=\frac{\sqrt{a^*}\eps _k}{\lam}u_k\Big(\frac{\eps _k}{\lam}x+x_k\Big) ,\,\ \text{where}\,\ \eps _k:=(a^*-a_k)^{ \frac{1}{2+p}}>0,
\]
so that   $\bar u_k(x)\to w(x)$ uniformly in $\R^2$ as $k\to\infty$. We then derive from (\ref{GNineq}) that
\begin{align}\label{3:eqa}
e(a_k) = E_{a_k}(u_k) & = \frac {\lam ^2} {a^*\eps_k^2} \Big[  \inte |\nabla
\bar{u}_k (x)|^2 dx - \frac {1} 2 \inte \bar{u}_k ^4(x) dx
\Big]+ \frac{\lam ^2\eps_k^p}{2(a^*)^2} \inte \bar{u}^4_k (x) dx\nonumber
\\ & + \frac{1}{a^*}\inte V\big(\frac{\eps_k}{\lam} x+x_k ) \bar{u}_k^2 (x) dx\\
& \ge  \frac{\lam ^2\eps_k^p}{2(a^*)^2} \inte \bar{u}^4_k (x) dx+ \frac{1}{a^*}\Big(\frac{\eps_k}{\lam}\Big)^p\inte g\big(\frac{\eps_k}{\lam} x+x_k )h\big(x+\frac{\lam x_k}{\eps_k } ) \bar{u}_k^2 (x) dx, \nonumber
\end{align}
which then implies from  (\ref{2:cona}) that $|\frac{\lam x_k}{\eps_k }|$ is bounded uniformly in $k$. Therefore, there exist  a subsequence (still denoted by $\{\frac{\lam x_k}{\eps_k }\}$) of $\{\frac{\lam x_k}{\eps_k }\}$ and $y^0\in\R^2$ such that
\[
\frac{\lam x_k}{\eps_k }\to y^0\,\ \text{as}\,\ k\to \infty.
\]
 Note that
\begin{equation}\label{strr} \arraycolsep=1.5pt\begin{array}{lll}
 && \liminf_{k \to \infty} \displaystyle \inte g\big(\frac{\eps_k}{\lam}x+x_k )h\big(x+\frac{\lam x_k}{\eps_k } )  \bar{u}_k^2 (x) dx \\[3mm]
 & \ge&  \liminf_{k \to \infty}  \displaystyle\int _{B_{\frac{1}{\sqrt {\eps_k}}}(0)} g\big(\frac{\eps_k}{\lam}x+x_k )h\big(x+\frac{\lam x_k}{\eps_k } )  \bar{u}_k^2 (x) dx\\[3mm] &= &\displaystyle g(0)\inte h (x+y^0 ) w ^2 (x) dx.
 \end{array}\end{equation}
Since $u_k$ gives the least energy of $e(a_k)$ and the assumption (\ref{1:HH}) implies that $y_0$ is essentially the unique global minimum point of $H(y)=\inte h (x+y ) w ^2 (x) dx$, we conclude from (\ref{3:eqa}) and (\ref{strr}) that $y^0=y_0 $, which thus implies that (\ref{2:conb}) holds, and the proof is therefore complete.
\qed

\subsection{Proof of local uniqueness}

Following Proposition \ref{2:prop}, this subsection is focussed on the proof of Theorem \ref{1:thmA}, and in the whole subsection we always assume that $V(x)=h(x)\in C^2 (\R^2)$ is homogeneous of degree $p\ge 2$ and satisfies (\ref{1:H}) and $\lim_{|x|\to\infty} h(x) = \infty$. Our proof is stimulated by \cite{Cao,Deng,Grossi}.
We first define the linearized operator $\mathcal{L}$ by
\[
\mathcal{L}:=-\Delta +(1-3w^2)\ \  \mbox{in}\ \, \R^2,
\]
where $w=w(|x|)>0$ is the unique positive  solution of (\ref{Kwong}) and $w$ satisfies the exponential decay (\ref{1:exp}). Recall from \cite{K,NT} that
\begin{equation}
ker (\mathcal{L})=span \Big\{\frac{\partial w}{\partial x_1},\frac{\partial w}{\partial x_2}\Big\}.
\label{2:linearized}
\end{equation}
For any positive minimizer $u_k$  of $e(a_k)$, where $a_k\nearrow a^*$ as $k\to\infty$, one can note that $u_k$ solves the Euler-Lagrange equation
\begin{equation}
-\Delta u_k(x) +V(x)u_k(x) =\mu_k u_k(x) +a_ku ^3 _k(x)\ \, \text{in\  $\R^2$},
\label{2:cond}
\end{equation}
where $\mu _k\in \R$ is a suitable Lagrange multiplier and satisfies
\begin{equation}\label{2:cone}
\mu_k = e(a_k) - \frac{a_k}{2} \inte u^4_k(x) dx.
\end{equation}
Moreover, under the more general assumption  (\ref{2:V1}), one can derive from (\ref{1:id}) and (\ref{2:cona}) that $u_k$ satisfies
\begin{equation}\label{2:conf}
 \inte u^4_k(x) dx=(a^*-a_k)^{-\frac{2}{2+p}}\Big[\frac{2\lam ^2}{a^*}+o(1)\Big] \,\ \mbox{as} \,\ k\to\infty .
\end{equation}
It then follows from (\ref{1:id}), (\ref{2:cone}) and (\ref{2:conf}) that $\mu_k$ satisfies
\begin{equation}\label{2:conmu}
\frac{\mu _k\eps _k^2}{\lam ^2} \to -1\, \ \mbox{as} \,\ k \to +\infty,
\end{equation}
where we denote
\[
\eps _k:=(a^*-a_k)^{ \frac{1}{2+p}}>0.
\]

Set
\[\bar u_k(x):=\frac{\sqrt{a^*}\eps _k}{\lam}u_k\Big(\frac{\eps _k}{\lam}x+x_k\Big) ,\]
so that Proposition \ref{2:prop} gives $\bar u_k(x)\to w(x)$ uniformly in $\R^2$ as $k\to\infty$. Note from (\ref{2:cond}) that $ \bar u_k$ satisfies
\begin{equation}\label{naa:A}
-\Delta \bar u_k(x)+\Big(\frac{\eps _k }{\lam  }\Big)^2V\Big(\frac{\eps _k}{\lam}x+x_k\Big)\bar u_k(x)=\frac{\mu _k\eps _k^2}{\lam ^2}\bar u_k(x)+\frac{a_k}{a^*}\bar u^3_k(x)\,\ \mbox{in} \,\ \R^2.
\end{equation}
Moreover, by the exponential decay (\ref{2:conexp}), there exist $C_0>0$ and $R>0$ such that
\begin{equation}\label{na:A}
|\bar u_k(x)|\le C_0 e^{-\frac{|x|}{2}}\,\ \mbox{for}  \,\ |x|>R,
\end{equation}
which then implies that
\[ \Big | \Big(\frac{\eps _k }{\lam  }\Big)^2V\Big(\frac{\eps _k}{\lam}x+x_k\Big)\bar u_k(x) \Big| \le  CC_0 e^{-\frac{|x|}{4}}\,\ \mbox{for}  \,\ |x|>R, \]
if $V(x)$ satisfies (\ref{2:V1}) with $p\ge 2$. Therefore, under the assumption  (\ref{2:V1}), applying the local elliptic estimates (see (3.15) in \cite{GT}) to (\ref{naa:A}) yields that
\begin{equation}\label{na:B}
|\nabla \bar u_k(x)|\le Ce^{- \frac{|x|}{4}}
\,\ \mbox{as} \,\ |x|\to\infty ,
\end{equation}
where the estimates (\ref{2:conmu}) and (\ref{na:A}) are also used. In the following, we shall follow Proposition \ref{2:prop} and (\ref{na:B}) to derive Theorem \ref{1:thmA} on the local uniqueness of positive minimizers as $a\nearrow  a^*$.

\vskip 0.1truein

\noindent\textbf{Proof of Theorem \ref{1:thmA}.}
Suppose that there exist two different positive minimizers $u_{1,k}$ and $u_{2,k}$ of $e(a_k)$ with $a_k\nearrow a^*$ as $k\to\infty$. Let $x_{1,k}$ and $x_{2,k}$ be the unique local maximum point of $u_{1,k}$ and $u_{2,k}$, respectively. Following (\ref{2:cond}), $u_{i,k}$ then solves the Euler-Lagrange equation
\begin{equation}
-\Delta u_{i,k}(x) +h(x)u_{i,k}(x) =\mu_{i,k} u_{i,k}(x) +a_ku ^3 _{i,k}(x)\quad \text{in\ \ $\R^2$}, \ \ i=1,2,
\label{3d:cond}
\end{equation}
where $V(x)=h(x)$ and $\mu _{i,k}\in \R$ is a suitable Lagrange multiplier.
Define
\begin{equation}\label{3d:27A}
\bar u_{i,k}(x):=\frac{\sqrt{a^*}\eps _k}{\lam}u_{i,k} ( x ), \ \ \mbox{where}\ \ i=1,2.
\end{equation}
Proposition \ref{2:prop} then implies that $\bar u_{i,k}\big(\frac{\eps _k}{\lam}x+x_{2,k}\big)\to w(x)$ uniformly in $\R^2$, and $\bar u_{i,k}$ satisfies the equation
\begin{equation}
-\eps ^2 _k\Delta \bar u_{i,k}(x) +\eps ^2 _kh(x)\bar u_{i,k}(x) =\mu_{i,k}\eps ^2 _k\bar  u_{i,k}(x) +\frac{\lam ^2a_k}{a^*}\bar u ^3 _{i,k}(x)\quad \text{in\ \ $\R^2$}, \ \ i=1,2.
\label{5.2:0}
\end{equation}
Because $u_{1,k}\not\equiv u_{2,k}$, we consider
\[
\bar \xi_k(x)=\frac{ u_{2,k}(x)-  u_{1,k}(x)}{\|  u_{2,k}-  u_{1,k}\|_{L^\infty(\R^2)}}=\frac{\bar u_{2,k}(x)-\bar u_{1,k}(x)}{\|\bar u_{2,k}-\bar u_{1,k}\|_{L^\infty(\R^2)}}.
\]
Then $\bar \xi_k$ satisfies the equation
\begin{equation}
-\eps ^2_k\Delta \bar \xi_k +\bar C_{k}(x)\bar \xi_k =\bar g_k(x)\quad \text{in\,\, $\R^2$},
\label{5.2:1}
\end{equation}
where the coefficient $\bar C_k(x)$ satisfies
\begin{equation}
\bar C_k(x):=-\mu _{1,k}\eps ^2_k- \frac{\lam ^2a_k}{a^*} \big(\bar u_{2,k}^2+\bar u_{2,k}\bar u_{1,k}+\bar u_{1,k}^2\big)+ \eps ^2_k  h(x),
\label{5.2:2}
\end{equation}
and the nonhomogeneous term $\bar g_k(x)$ satisfies
\begin{equation}
 \arraycolsep=1.5pt\begin{array}{lll}
\bar g_k(x):=\displaystyle \frac{\eps ^2_k\bar u_{2,k}(\mu _{2,k}-\mu _{1,k})}{\|\bar u_{2,k}-\bar u_{1,k}\|_{L^\infty(\R^2)}}&=&-\displaystyle  \frac{\lam ^4a_k\bar u_{2,k}}{2(a^*)^2\eps ^2_k} \inte \frac{\bar u_{2,k}^4-\bar u_{1,k}^4 }{\|\bar u_{2,k}-\bar u_{1,k}\|_{L^\infty(\R^2)}}dx\\[4mm]
&=&-\displaystyle  \frac{\lam ^4a_k\bar u_{2,k}}{2(a^*)^2\eps ^2_k}\inte \bar \xi_k\big(\bar u_{2,k}^2+\bar u_{1,k}^2\big)\big(\bar u_{2,k}+\bar u_{1,k}\big)dx,
\end{array}
\label{5.2:3}
 \end{equation}
due to the relation (\ref{2:cone}).

Motivated by \cite{Cao}, we first claim that for any $x_0\in\R^2$, there exists a small constant $\delta >0$  such that
\begin{equation}
    \int_{\partial B_\delta (x_0)} \Big[ \eps ^2_k |\nabla \bar \xi_k|^2+ \frac{\lam ^2}{2} |\bar\xi_k|^2+ \eps ^2_k  h(x)|\bar\xi_k|^2\Big]dS=O( \eps ^2_k)\quad \text{as}\ k\to\infty.
\label{5.2:6}
\end{equation}
To prove the above claim, multiplying (\ref{5.2:1}) by $\bar \xi_k$ and integrating over $\R^2$, we obtain that
\[\arraycolsep=1.5pt\begin{array}{lll}
&&\displaystyle \eps ^2_k\inte |\nabla \bar \xi_k|^2 -\mu_{i,k}\eps ^2 _k\inte  |\bar\xi_k|^2+\eps ^2_k\inte h(x)|\bar\xi_k|^2 \\[4mm]
&=&\displaystyle \frac{\lam ^2a_k}{a^*}\inte \big(\bar u_{2,k}^2+\bar u_{2,k}\bar u_{1,k}+\bar u_{1,k}^2\big)|\bar\xi_k|^2\\[4mm]
&&-
\displaystyle  \frac{\lam ^4a_k}{2(a^*)^2\eps ^2_k}\inte \bar u_{2,k}\bar\xi_k\inte \bar \xi_k\big(\bar u_{2,k}^2+\bar u_{1,k}^2\big)\big(\bar u_{2,k}+\bar u_{1,k}\big)
\\[4mm]
&\le &\displaystyle \frac{\lam ^2a_k}{a^*}\inte \big(\bar u_{2,k}^2+\bar u_{2,k}\bar u_{1,k}+\bar u_{1,k}^2\big)+\displaystyle  \frac{\lam ^4a_k}{2(a^*)^2\eps ^2_k}\inte \bar u_{2,k} \inte  \big(\bar u_{2,k}^2+\bar u_{1,k}^2\big)\big(\bar u_{2,k}+\bar u_{1,k}\big)\\[4mm]
&\le& C\eps ^2_k\,\ \mbox{as} \,\ k\to\infty,
\end{array}\]
since $|\bar\xi_k|$ and $\bar u_{i,k}\big(\frac{\eps _k}{\lam}x+x_{2,k}\big)$ are bounded uniformly in $k$, and $\bar u_{i,k}\big(\frac{\eps _k}{\lam}x+x_{2,k}\big)$ decays exponentially as $|x|\to\infty$, $i=1,\, 2$. This implies that there exists a constant $C_1>0$ such that
\begin{equation}
     I:=\eps ^2_k\inte |\nabla \bar \xi_k|^2+\frac{\lam ^2}{2} \inte |\bar\xi_k|^2+ \eps ^2_k\inte h(x)|\bar\xi_k|^2<C_1\eps ^2_k \quad \text{as}\ \, k\to\infty.
\label{5.2:5}
\end{equation}
Applying Lemma 4.5 in \cite{Cao}, we then conclude that for any $x_0\in\R^2$, there exist a small constant $\delta >0$ and $C_2>0$  such that
\[
    \int_{\partial B_\delta (x_0)} \Big[ \eps ^2_k |\nabla \bar \xi_k|^2+ \frac{\lam ^2}{2}  |\bar\xi_k|^2+ \eps ^2_k  h(x)|\bar\xi_k|^2\Big]dS\le C_2I\le C_1C_2\eps ^2_k\,\ \mbox{as} \,\ k\to\infty,
\]
which therefore implies the claim (\ref{5.2:6}).

We next define
\begin{equation}
 \xi _k(x)=\bar \xi_k\big(\frac{\eps _k}{\lam}x+x_{2,k}\big), \ \ k=1,2,\cdots ,
\label{5.2:xi}
\end{equation}
and
\[
\tilde{u}_{i,k}(x):=\frac{\sqrt{a^*}\eps _k}{\lam}u_{i,k} \big(\frac{\eps _k}{\lam}x+x_{2,k}\big), \ \ \mbox{where}\ \ i=1,2,
\]
so that   $\tilde{u}_{i,k}(x)\to w(x)$ uniformly in $\R^2$ as $k\to\infty$ in view of Proposition \ref{2:prop}.
Under the non-degeneracy assumption (\ref{1:H}), we shall carry out the proof of Theorem \ref{1:thmA} by deriving a contradiction through the following three steps.

\vskip 0.1truein

\noindent{\em  Step 1.} There exist a subsequence $\{a_k\}$ and some constants $b_0$, $b_1$ and $b_2$ such that $\xi_k(x)\to \xi_0(x)$ in $C_{loc}(\R^2)$ as $k\to\infty$, where
\begin{equation}
\xi_0(x)=b_0\big(w+x\cdot \nabla w\big)+\sum ^2_{i=1}b_i\frac{\partial w}{\partial x_i}.
\label{3d:5a}
\end{equation}

Note that $\xi_k$ satisfies
\begin{equation}
-\Delta \xi_k +C_{k}(x)\xi_k =g_k(x)\quad \text{in\,\, $\R^2$},
\label{3d:2}
\end{equation}
where the coefficient $C_k(x)$ satisfies
\begin{equation}\arraycolsep=1.5pt\begin{array}{lll}
C_k(x)&:=&-\Big(1-\displaystyle\frac{\eps_k^{2+p}}{a^*}\Big)\Big[\tilde{u}_{2,k}^2(x)+\tilde{u}_{2,k}(x)\tilde{u}_{1,k}(x)+
\tilde{u}_{1,k}^2(x)\Big]\\[4mm]
&&-\displaystyle\frac{\eps ^2_k}{\lam ^2}\mu _{1,k}+\displaystyle\frac{\eps ^2_k}{\lam ^2}h\big(\frac{\eps_kx}{\lam}+x_{2,k}\big),
\end{array}\label{3d:3}
\end{equation}
and the nonhomogeneous term $g_k(x)$ satisfies
\begin{equation}
 \arraycolsep=1.5pt\begin{array}{lll}
g_k(x):=\displaystyle\frac{\tilde{u}_{2,k}}{\lam ^2}\frac{\eps ^2_k(\mu _{2,k}-\mu _{1,k})}{\|\tilde{u}_{2,k}-\tilde{u}_{1,k}\|_{L^\infty}}&=&-\displaystyle \frac{\tilde{u}_{2,k}}{\lam ^2}\frac{a_k\eps ^2_k}{2} \inte \frac{ u_{2,k}^4-u_{1,k}^4 }{\|\tilde{u}_{2,k}-\tilde{u}_{1,k}\|_{L^\infty}}dx\\[4mm]
&=&-\displaystyle \frac{a_k\tilde{u}_{2,k}}{2(a^*)^2}\inte \xi_k\big(\tilde{u}_{2,k}^2+\tilde{u}_{1,k}^2\big)\big(\tilde{u}_{2,k}+\tilde{u}_{1,k}\big)dx.
\end{array}
\label{3d:4}
 \end{equation}
Here we have used (\ref{2:cone}) and (\ref{5.2:0}). Since $\|\xi _k\|_{L^\infty(\R^2)}\le 1$, the standard elliptic regularity then implies (cf. \cite{GT}) that $\|\xi _k\|_{C^{1,\alpha }_{loc}(\R^2)}\le C$ for some $\alp \in (0,1)$, where the constant $C>0$ is independent of $k$. Therefore, there exist  a subsequence $\{a_k\}$ and a function $\xi _0=\xi _0(x)$ such that $\xi_k(x)\to \xi_0(x)$ in $C_{loc}(\R^2)$ as $k\to\infty$.
Applying  Proposition \ref{2:prop}, direct calculations yield from (\ref{2:cone}) and (\ref{2:conf}) that
\[
C_k(x)\to 1-3w^2(x)\ \ \text{uniformly\ on\ $\R^2$}\,\ \mbox{as} \,\ k\to\infty,
\]
and
\[
g_k(x)\to -\frac{2w(x)}{a^*}\inte w^3\xi _0\ \ \text{uniformly\ on\ $\R^2$}\,\ \mbox{as} \,\ k\to\infty.
\]
This implies from (\ref{3d:2}) that $\xi_0$ solves
\begin{equation}
\mathcal{L}\xi_0=-\Delta \xi_0+(1-3w^2)\xi_0=\Big(-\frac{2}{a^*}\inte w^3\xi _0\Big)w\ \ \mbox{in} \ \ \R^2.
\label{3d:5}
\end{equation}
Since $\mathcal L(w+ x\cdot \nabla w)=-2 w$,   we then conclude from (\ref{2:linearized}) and (\ref{3d:5}) that (\ref{3d:5a}) holds for some constants $b_0$, $b_1$ and $b_2$.

\vskip 0.1truein

\noindent{\em  Step 2.} The constants $b_0=b_1=b_2=0$ in (\ref{3d:5a}).

We first derive the following Pohozaev-type identity
\begin{equation}
b_0 \displaystyle\inte \frac{\partial h(x+y_0)}{\partial x_j}\big(x\cdot \nabla w^2\big) -\displaystyle\sum ^2_{i=1}b_i\inte \frac{\partial ^2 h(x+y_0)}{\partial x_j\partial x_i}w^2=0,\quad j=1,\,2.
\label{5.2:AA}
\end{equation}
Multiplying (\ref{5.2:0}) by $\frac{\partial \bar u_{i,k}}{\partial  x_j}$, where $i,j=1,2$, and integrating over $B_\delta (x_{2,k})$, where $\delta >0$ is small and given by (\ref{5.2:6}), we calculate that
\begin{equation}\arraycolsep=1.5pt\begin{array}{lll}
&&-\eps _k^2\displaystyle\intB\frac{\partial \bar u_{i,k}}{\partial  x_j}\Delta \bar u_{i,k}+\eps _k^2\displaystyle\intB h(x)\frac{\partial \bar u_{i,k}}{\partial  x_j} \bar u_{i,k}\\[4mm]
&=&\mu_{i,k}\eps _k^2\displaystyle\intB \frac{\partial \bar u_{i,k}}{\partial  x_j} \bar u_{i,k}+\displaystyle\frac{\lam ^2a_k}{a^*}\intB \frac{\partial \bar u_{i,k}}{\partial  x_j} \bar u_{i,k}^3\\[4mm]
&=&\displaystyle\frac{1}{2}\mu_{i,k}\eps _k^2\intPB \bar u_{i,k}^2\nu _jdS+\displaystyle\frac{\lam ^2a_k}{4a^*}\intPB \bar u_{i,k}^4\nu _jdS,
\end{array}\label{5.2:7}
\end{equation}
where $\nu =(\nu _1,\nu _2)$ denotes the outward unit normal of $\partial B_\delta (x_{2,k})$.
Note that
\[\arraycolsep=1.5pt\begin{array}{lll}
 &&-\eps _k^2\displaystyle\intB\frac{\partial \bar u_{i,k}}{\partial  x_j}\Delta \bar u_{i,k}\\[4mm]&=&-\eps _k^2\displaystyle\intPB\frac{\partial \bar u_{i,k}}{\partial  x_j}\frac{\partial \bar u_{i,k}}{\partial  \nu}dS+\eps _k^2\displaystyle\intB\nabla \bar u_{i,k}\cdot\nabla\frac{\partial \bar u_{i,k}}{\partial  x_j}\\[4mm]
 &=&-\eps _k^2\displaystyle\intPB\frac{\partial \bar u_{i,k}}{\partial  x_j}\frac{\partial \bar u_{i,k}}{\partial  \nu}dS+\displaystyle\frac{1}{2}\eps _k^2\intPB |\nabla \bar u_{i,k}|^2\nu _jdS,
\end{array}
\]
and
\[
 \eps _k^2\displaystyle\intB h(x)\frac{\partial \bar u_{i,k}}{\partial  x_j} \bar u_{i,k}= \frac{\eps _k^2}{2}\intPB h(x)\bar u_{i,k}^2\nu _jdS-\frac{\eps _k^2}{2}\intB \frac{\partial h(x)}{\partial  x_j}\bar u_{i,k}^2.
\]
We then derive from (\ref{5.2:7}) that
\begin{equation}\arraycolsep=1.5pt\begin{array}{lll}
&&\displaystyle\eps _k^2 \intB \frac{\partial h(x)}{\partial  x_j}\bar u_{i,k}^2\\[4mm]
&=&-2\eps _k^2\displaystyle\intPB\frac{\partial \bar u_{i,k}}{\partial  x_j}\frac{\partial \bar u_{i,k}}{\partial  \nu}dS+\displaystyle \eps _k^2\intPB |\nabla \bar u_{i,k}|^2\nu _jdS \\[4mm]
&& +\displaystyle\eps _k^2 \intPB h(x)\bar u_{i,k}^2\nu _jdS-\displaystyle \mu_{i,k}\eps _k^2\intPB \bar u_{i,k}^2\nu _jdS\\[4mm]
&&-\displaystyle\frac{\lam ^2a_k}{2a^*}\intPB \bar u_{i,k}^4\nu _jdS.
\end{array}\label{5.2:8}
\end{equation}
Following (\ref{5.2:8}), we thus have
\begin{equation}\arraycolsep=1.5pt\begin{array}{lll}
&&\displaystyle\eps _k^2 \intB \frac{\partial h(x)}{\partial  x_j}\big(\bar u_{2,k}+\bar u_{1,k}\big) \bar \xi _k dx\\[4mm]
&=&-2\displaystyle \eps _k^2\intPB \Big[\frac{\partial \bar u_{2,k}}{\partial  x_j}\frac{\partial \bar \xi_k}{\partial  \nu}+\frac{\partial \bar \xi_k}{\partial  x_j}\frac{\partial \bar u_{1,k}}{\partial  \nu}\Big]dS\\[4mm]
&&+\eps _k^2\displaystyle\intPB\nabla \bar \xi_k \cdot\nabla \big(\bar u_{2,k}+\bar u_{1,k}\big)\nu _jdS \\[4mm]
&& +\displaystyle \eps _k^2 \intPB h(x)\big(\bar u_{2,k}+\bar u_{1,k}\big) \bar \xi _k \nu _jdS-\displaystyle \mu_{1,k}\eps _k^2\intPB \big(\bar u_{2,k}+\bar u_{1,k}\big) \bar \xi _k\nu _jdS \\[4mm]
&& -\displaystyle\frac{\lam ^2a_k}{2a^*}\intPB \big(\bar u^2_{2,k}+\bar u^2_{1,k}\big) \big(\bar u_{2,k}+\bar u_{1,k}\big)\bar \xi _k\nu _jdS\\[4mm]
&&-\displaystyle\frac{\big(\mu_{2,k}-\mu_{1,k}\big)\eps _k^2}{\|\bar u_{2,k}-\bar u_{1,k}\|_{L^\infty}} \displaystyle\intPB \bar u_{2,k}^2\nu _jdS.
\end{array}\label{5.2:9}
\end{equation}

We now estimate the right hand side of (\ref{5.2:9}) as follows. Applying (\ref{5.2:6}), if $\delta >0$ is small, we then deduce that
\begin{equation}\arraycolsep=1.5pt\begin{array}{lll}
&&\displaystyle \eps _k^2\intPB \Big|\frac{\partial \bar u_{2,k}}{\partial  x_j}\frac{\partial \bar \xi_k}{\partial  \nu}\Big|dS\\[4mm]
&\le &\displaystyle \eps _k\Big(\intPB \Big|\frac{\partial \bar u_{2,k}}{\partial  x_j}\Big|^2dS\Big)^{\frac{1}{2}}\Big(\eps _k^2\intPB \Big|\frac{\partial \bar \xi_k}{\partial  \nu}\Big|^2dS\Big)^{\frac{1}{2}}\le C\eps _k^2e^{-\frac{C\delta}{\eps_k}}\,\ \mbox{as} \,\ k\to\infty,
\end{array}\label{5.2:9a}
\end{equation}
due to the fact that $\nabla \bar u_{2,k}\big(\frac{\eps _k}{\lam}x+x_{2,k}\big)$ satisfies the exponential decay (\ref{na:B}), where $C>0$ is independent of $k$. Similarly, we have
\[
\eps _k^2\intPB \Big|\frac{\partial \bar \xi_k}{\partial  x_j}\frac{\partial \bar u_{1,k}}{\partial  \nu} \Big|dS\le C\eps _k^2e^{-\frac{C\delta}{\eps_k}}\,\ \mbox{as} \,\ k\to\infty,
\]
and
\[
\eps _k^2\Big|\displaystyle\intPB \nabla \bar \xi_k \cdot\nabla \big(\bar u_{2,k}+\bar u_{1,k}\big)\nu _jdS\Big|\le C\eps _k^2e^{-\frac{C\delta}{\eps_k}}\,\ \mbox{as} \,\ k\to\infty,
\]
On the other hand, since both $|\bar\xi_k| $ and $|\big(\mu_{2,k}-\mu_{1,k}\big)\eps _k^2|$ are bounded uniformly in $k$, we also get from  (\ref{na:B}) that
 \begin{equation}\arraycolsep=1.5pt\begin{array}{lll}
&&\Big|\displaystyle \eps _k^2 \intPB h(x)\big(\bar u_{2,k}+\bar u_{1,k}\big) \bar \xi _k \nu _jdS-\displaystyle \mu_{1,k}\eps _k^2\intPB \big(\bar u_{2,k}+\bar u_{1,k}\big) \bar \xi _k\nu _jdS \\[4mm]
&& -\displaystyle\frac{\lam ^2a_k}{2a^*}\intPB \big(\bar u^2_{2,k}+\bar u^2_{1,k}\big) \big(\bar u_{2,k}+\bar u_{1,k}\big)\bar \xi _k\nu _jdS\\[4mm]
&&-\displaystyle\frac{\big(\mu_{2,k}-\mu_{1,k}\big)\eps _k^2}{\|\bar u_{2,k}-\bar u_{1,k}\|_{L^\infty}} \displaystyle\intPB \bar u_{2,k}^2\nu _jdS\Big|\\[4mm]
&=&o(e^{-\frac{C\delta}{\eps_k}})\,\ \mbox{as} \,\ k\to\infty,
\end{array}
\label{5.2:9b}
\end{equation}
due to the fact that  (\ref{5.2:3}) gives
\begin{equation}
 \displaystyle\frac{\big|\mu_{2,k}-\mu_{1,k}\big|\eps _k^2}{\|\bar u_{2,k}-\bar u_{1,k}\|_{L^\infty}}
 \le  \displaystyle\frac{\lam^4 a_k}{2(a^*)^2\eps _k^2}\inte \big(\bar u_{2,k}^2+\bar u_{1,k}^2\big)\big(\bar u_{2,k}+\bar u_{1,k}\big)|\bar \xi_k |
 \le  M ,
\label{5.2:9F}
\end{equation}
where the constants $M>0$ is independent of $k$.
Because $h(x)$ is homogeneous of degree $p$,  it then follows from (\ref{5.2:9}) that for small $\delta >0$,
 \begin{equation}\label{5.2:10}\arraycolsep=1.5pt\begin{array}{lll}
o(e^{-\frac{C\delta}{\eps_k}})&=&\displaystyle\eps _k^2 \intB \frac{\partial h(x)}{\partial  x_j}\big[\bar u_{2,k}(x)+\bar u_{1,k}(x)\big] \bar \xi _k(x)dx\\[4mm]
&=&\displaystyle\frac{\eps _k^3}{\lam } \int_{B_{\frac{\lam\delta}{\eps_k}}(0)} \frac{\partial}{\partial  y_j} h\big(\frac{\eps_k}{\lam}y+x_{2,k}\big)\bar \xi _k \big(\frac{\eps_k}{\lam}y+x_{2,k}\big)\\[4mm]
&&\qquad\qquad\,\ \cdot \displaystyle\Big[\bar u_{2,k}\big(\frac{\eps_k}{\lam}y+x_{2,k}\big)+\bar u_{1,k} \big(\frac{\eps_k}{\lam}y+x_{2,k}\big)\Big] dy\\[4mm]
&=&\displaystyle\frac{\eps _k^{p+3}}{\lam ^{p+1}} \Big[\int_{B_{\frac{\lam\delta}{\eps_k}}(0)} \frac{\partial}{\partial  y_j} h\big(y+\frac{\lam x_{2,k}}{\eps_k}\big)\bar \xi _k \big(\frac{\eps_k}{\lam}y+x_{2,k}\big)\\[4mm]
&&\qquad \ \cdot\displaystyle\Big[\bar u_{2,k}\big(\frac{\eps_k}{\lam}y+x_{2,k}\big)+\bar u_{1,k} \big(\frac{\eps_k}{\lam}y+x_{1,k}\big)\Big] dy+o (1)\Big]
\end{array}\end{equation}
as $k\to\infty$. Applying (\ref{1:H}), we thus derive from (\ref{2:conb}), (\ref{3d:5a}) and (\ref{5.2:10}) that
\[ \arraycolsep=1.5pt\begin{array}{lll}
0=2\displaystyle\inte \frac{\partial h(x+y_0)}{\partial x_j}w\,\xi_0&=&\displaystyle2\inte \frac{\partial h(x+y_0)}{\partial x_j}w\Big[b_0\big(w+x\cdot \nabla w\big)+\sum ^2_{i=1}b_i\frac{\partial w}{\partial x_i}\Big]\\[4mm]
&=&b_0 \displaystyle\inte \frac{\partial h(x+y_0)}{\partial x_j}\big(x\cdot \nabla w^2\big) -\displaystyle\sum ^2_{i=1}b_i\inte \frac{\partial ^2 h(x+y_0)}{\partial x_j\partial x_i}w^2,
\end{array}\]
where $j=1,\,2$, which thus implies (\ref{5.2:AA}).


\vskip 0.1truein

We next derive $b_0=0$. Using the integration by parts,  we note that
\begin{equation}\arraycolsep=1.5pt\begin{array}{lll}
&&-\displaystyle\eps _k^2 \intB \big[(x-x_{2,k})\cdot \nabla \bar u_{i,k}\big] \Delta \bar u_{i,k} \\[4mm]
&=&- \eps _k^2\displaystyle \intPB \frac{\partial\bar u_{i,k} }{\partial \nu }\big[(x-x_{2,k})\cdot \nabla \bar u_{i,k}\big]
+\displaystyle\eps _k^2 \intB \nabla \bar u_{i,k}\nabla \big[(x-x_{2,k})\cdot \nabla \bar u_{i,k}\big] \\[4mm]
&=& - \eps _k^2\displaystyle \intPB \frac{\partial\bar u_{i,k} }{\partial \nu }\big[(x-x_{2,k})\cdot \nabla \bar u_{i,k}\big]
+\displaystyle \frac{\eps _k^2}{2}\intPB \big[(x-x_{2,k})\cdot \nu \big]|\nabla \bar u_{i,k}|^2.
\end{array}\label{5.3:3}
\end{equation}
Multiplying (\ref{5.2:0}) by $ (x-x_{2,k})\cdot \nabla \bar u_{i,k} $, where $i=1,2$, and integrating over $B_\delta (x_{2,k})$, where $\delta >0$ is small as before, we deduce that for $i=1,2,$
\begin{equation}\arraycolsep=1.5pt\begin{array}{lll}
&&-\displaystyle\eps _k^2 \intB \big[(x-x_{2,k})\cdot \nabla \bar u_{i,k}\big] \Delta \bar u_{i,k} \\[4mm]
&=& \displaystyle\eps _k^2 \intB \big[\mu_{i,k}-h(x)\big] \bar u_{i,k}\big[(x-x_{2,k})\cdot \nabla \bar u_{i,k}\big]\\[4mm]
&& +
\displaystyle \frac{\lam ^2 a_k}{a^*} \intB \bar u_{i,k}^3\big[(x-x_{2,k})\cdot \nabla \bar u_{i,k}\big]\\[4mm]
&=& -\displaystyle\frac{\eps _k^2}{2} \intB\bar u_{i,k}^2\Big\{2\big[\mu_{i,k}-h(x)\big]-(x-x_{2,k})\cdot \nabla h(x)\Big\}\\[4mm]
&&+\displaystyle\frac{\eps _k^2}{2} \intPB\bar u_{i,k}^2\big[\mu_{i,k}-h(x)\big](x-x_{2,k})\nu dS\\[4mm]
&&-\displaystyle \frac{\lam ^2 a_k}{2a^*} \intB \bar u_{i,k}^4+\displaystyle \frac{\lam ^2 a_k}{4a^*} \intPB \bar u_{i,k}^4(x-x_{2,k})\nu dS\\[4mm]
&=& -\mu_{i,k}\displaystyle\eps _k^2 \inte \bar u_{i,k}^2+\displaystyle\frac{2+p}{2}\eps _k^2 \inte h(x)\bar u_{i,k}^2-\displaystyle \frac{\lam ^2 a_k}{2a^*} \inte \bar u_{i,k}^4+I_i,
\end{array}\label{5.3:1}
\end{equation}
where the lower order term $I_i$ satisfies
\begin{equation}\arraycolsep=1.5pt\begin{array}{lll}
I_i&=& \mu_{i,k}\displaystyle\eps _k^2 \int _{\R^2\backslash B_\delta (x_{2,k})} \bar u_{i,k}^2-\displaystyle\frac{2+p}{2}\eps _k^2 \int _{\R^2\backslash B_\delta (x_{2,k})} h(x)\bar u_{i,k}^2\\[4mm]
&&+\displaystyle \frac{\lam ^2 a_k}{2a^*} \int _{\R^2\backslash B_\delta (x_{2,k})} \bar u_{i,k}^4-\displaystyle\frac{1}{2} \eps _k^2\intB \bar u_{i,k}^2\big[x_{2,k}\cdot \nabla h(x)\big]\\[4mm]
&&+\displaystyle\frac{\eps _k^2}{2} \intPB\bar u_{i,k}^2\big[\mu_{i,k}-h(x)\big](x-x_{2,k})\nu dS\\[4mm]
&&+\displaystyle \frac{\lam ^2 a_k}{4a^*} \intPB \bar u_{i,k}^4(x-x_{2,k})\nu dS,\,\ i=1,2.
\end{array}\label{5.3:2}
\end{equation}
Since it follows from (\ref{2:cone}) that
\[
-\mu_{i,k}\displaystyle\eps _k^2 \inte \bar u_{i,k}^2-\displaystyle \frac{\lam ^2 a_k}{2a^*} \inte \bar u_{i,k}^4=-\frac{a^*\eps_k^4}{\lam ^2}\Big[\mu_{i,k}+\frac{a_k}{2}\inte  u_{i,k}^4\Big]=-\frac{a^*\eps_k^4}{\lam ^2}e(a_k),
\]
we reduce from (\ref{5.3:3})--(\ref{5.3:2}) that
\[\arraycolsep=1.5pt\begin{array}{lll}
\displaystyle\frac{a^*\eps_k^4}{\lam ^2}e(a_k)-\displaystyle\frac{2+p}{2}\eps _k^2 \inte h(x)\bar u_{i,k}^2&=&I_i+\eps _k^2\displaystyle \intPB \frac{\partial\bar u_{i,k} }{\partial \nu }\big[(x-x_{2,k})\cdot \nabla \bar u_{i,k}\big]\\[4mm]
&&-\displaystyle \frac{\eps _k^2}{2}\intPB \big[(x-x_{2,k})\cdot \nu \big]|\nabla \bar u_{i,k}|^2,\,\ i=1,2,
\end{array}\]
which implies that
\begin{equation}
-\displaystyle\frac{2+p}{2}\eps _k^2 \inte h(x)\big[\bar u_{2,k}+\bar u_{1,k}\big]\bar \xi_k=T_k.
 \label{5.3:4}
\end{equation}
Here the term $T_k$ satisfies that for small $\delta >0$,
\begin{equation}\arraycolsep=1.5pt\begin{array}{lll}
T_k&=&\displaystyle\frac{I_2-I_1}{\|\bar u_{2,k}-\bar u_{1,k}\|_{L^\infty}}-\displaystyle \frac{\eps _k^2}{2}\intPB \big[(x-x_{2,k})\cdot \nu \big]\big(\nabla \bar u_{2,k}+\nabla \bar u_{1,k}\big)\nabla \bar \xi_k\\[4mm]
&&+\displaystyle  \eps _k^2 \intPB \Big\{\big[(x-x_{2,k})\cdot \nabla \bar u_{2,k}\big]\big(\nu \cdot \nabla\bar \xi_k\big)+\big(\nu \cdot \nabla \bar u_{1,k}\big)\big[(x-x_{2,k})\cdot \nabla \bar \xi_k\big]\Big\}
\\[4mm]
&=&\displaystyle\frac{I_2-I_1}{\|\bar u_{2,k}-\bar u_{1,k}\|_{L^\infty}}+o(e^{-\frac{C\delta}{\eps_k}})\,\ \mbox{as} \,\ k\to\infty,
\end{array} \label{5.3:5}\end{equation}
due to (\ref{na:B}) and (\ref{5.2:6}), where the second equality follows by applying the argument of estimating  (\ref{5.2:9a}).

Using the arguments of estimating (\ref{5.2:9a}) and (\ref{5.2:9b}), along with the exponential decay of $\bar u_{i,k}$, we also derive that for small $\delta >0$,
\begin{equation}\arraycolsep=1.5pt\begin{array}{lll}
&&\displaystyle\frac{I_2-I_1}{\|\bar u_{2,k}-\bar u_{1,k}\|_{L^\infty}}\\[4mm]
&=&\mu_{2,k}\displaystyle\eps _k^2 \int _{\R^2\backslash B_\delta (x_{2,k})} \big(\bar u_{2,k}+\bar u_{1,k}\big)\bar \xi_k-\displaystyle\frac{2+p}{2}\eps _k^2 \int _{\R^2\backslash B_\delta (x_{2,k})} h(x)\big(\bar u_{2,k}+\bar u_{1,k}\big)\bar \xi_k\\[4mm]
&&+\displaystyle \frac{\lam ^2 a_k}{2a^*} \int _{\R^2\backslash B_\delta (x_{2,k})} \big(\bar u_{2,k}^2+\bar u_{1,k}^2\big)\big(\bar u_{2,k}+\bar u_{1,k}\big)\bar \xi_k\\[4mm]
&&+\displaystyle\frac{\big(\mu_{2,k}-\mu_{1,k}\big)\eps _k^2}{\|\bar u_{2,k}-\bar u_{1,k}\|_{L^\infty}}  \displaystyle\int _{\R^2\backslash B_\delta (x_{2,k})}  \bar u_{1,k}^2-\displaystyle\frac{1}{2}\eps _k^2\intB \big[x_{2,k}\cdot \nabla h(x)\big]\big(\bar u_{2,k}+\bar u_{1,k}\big)\bar \xi_k\\[4mm]
&&+\displaystyle \frac{\lam ^2 a_k}{4a^*} \intPB   \big(\bar u_{2,k}^2+\bar u_{1,k}^2\big)\big(\bar u_{2,k}+\bar u_{1,k}\big)\bar \xi_k(x-x_{2,k})\nu dS\\[4mm]
&&-\displaystyle\frac{\eps _k^2}{2} \intPB \big(\bar u_{2,k}+\bar u_{1,k}\big)\bar \xi_kh(x)(x-x_{2,k})\nu dS\\[4mm]
&&+\displaystyle\frac{\mu_{2,k}\eps _k^2}{2} \intPB   \big(\bar u_{2,k}+\bar u_{1,k}\big)\bar \xi_k(x-x_{2,k})\nu dS\\[4mm]
&&
+\displaystyle\frac{\big(\mu_{2,k}-\mu_{1,k}\big)\eps _k^2}{2\|\bar u_{2,k}-\bar u_{1,k}\|_{L^\infty}} \intPB   \bar u_{1,k}^2 (x-x_{2,k})\nu dS\\[4mm]
&=&\displaystyle\frac{\big(\mu_{2,k}-\mu_{1,k}\big)\eps _k^2}{\|\bar u_{2,k}-\bar u_{1,k}\|_{L^\infty}} \Big[ \displaystyle\int _{\R^2\backslash B_\delta (x_{2,k})}  \bar u_{1,k}^2+\displaystyle\frac{1}{2} \intPB   \bar u_{1,k}^2 (x-x_{2,k})\nu dS\Big]\\[4mm]
&&-\displaystyle\frac{1}{2} \eps _k^2\intB \big[ x_{2,k} \cdot \nabla h(x)\big]\big(\bar u_{2,k}+\bar u_{1,k}\big)\bar \xi_k+o(e^{-\frac{C\delta}{\eps_k}})\,\ \mbox{as} \,\ k\to\infty.
\end{array} \label{5.3:10}
\end{equation}
Note from (\ref{5.2:9F}) that
\begin{equation}
\displaystyle\frac{\big(\mu_{2,k}-\mu_{1,k}\big)\eps _k^2}{\|\bar u_{2,k}-\bar u_{1,k}\|_{L^\infty}} \Big[ \displaystyle\int _{\R^2\backslash B_\delta (x_{2,k})}  \bar u_{1,k}^2+\displaystyle\frac{1}{2} \intPB   \bar u_{1,k}^2 (x-x_{2,k})\nu dS\Big] =O(e^{-\frac{C\delta}{\eps_k}})
 \label{5.3:10b}
\end{equation}
as $k\to\infty$, where the constant $C>0$ is independent of $k$. Moreover, we follow from the first identity of (\ref{5.2:10}) that
\begin{equation}\arraycolsep=1.5pt\begin{array}{lll}
&& \displaystyle\frac{1}{2} \eps _k^2\intB \big[ x_{2,k} \cdot \nabla h(x)\big]\big(\bar u_{2,k}+\bar u_{1,k}\big)\bar \xi_k \\[4mm]
&= & \displaystyle\frac{1}{2}  \eps _k^2\displaystyle\sum ^2_{i=1}x^i_{2,k}\intB \frac{\partial h(x)}{\partial  x_i}\big[\bar u_{2,k}(x)+\bar u_{1,k}(x)\big] \bar \xi _k(x)dx \\[4mm]
&=& o(e^{-\frac{C\delta}{\eps_k}})\,\ \mbox{as} \,\ k\to\infty,
\end{array}\label{5.3:10c}
\end{equation}
where we denote $x_{2,k}=(x^1_{2,k},x^2_{2,k})$.
Therefore, we deduce from (\ref{5.3:5})--(\ref{5.3:10c}) that
\[
T_k=o(\eps _k^{4+p}) \,\ \mbox{as} \,\ k\to\infty.
\]
Further, we obtain from (\ref{5.3:4}) that
\[\arraycolsep=1.5pt\begin{array}{lll}
o(\eps _k^{4+p})&=&-\displaystyle\frac{2+p}{2}\eps _k^2 \inte h(x)\Big[\bar u_{2,k}+\bar u_{1,k}\Big]\bar \xi_k\\[4mm]
&=&-\displaystyle\frac{2+p}{2\lam^2}\eps _k^4 \inte h\big(\frac{\eps _k}{\lam}x+x_{2,k}\big)\Big[\bar u_{2,k}\big(\frac{\eps _k}{\lam}x+x_{2,k}\big)+\bar u_{1,k}\big(\frac{\eps _k}{\lam}x+x_{1,k}\big)\Big]  \xi_k(x)dx\\[4mm]
&&-\displaystyle\frac{2+p}{2\lam^2}\eps _k^4 \inte h\big(\frac{\eps _k}{\lam}x+x_{2,k}\big)\Big[\bar u_{1,k}\big(\frac{\eps _k}{\lam}x+x_{2,k}\big)-\bar u_{1,k}\big(\frac{\eps _k}{\lam}x+x_{1,k}\big)\Big]  \xi_k(x)dx\\[4mm]
&=&-\displaystyle\frac{2+p}{2\lam^{2+p}}\eps _k^{4+p} \inte h\big(x+\frac{\lam x_{2,k}}{\eps _k}\big)\Big[\bar u_{2,k}\big(\frac{\eps _k}{\lam}x+x_{2,k}\big)+\bar u_{1,k}\big(\frac{\eps _k}{\lam}x+x_{1,k}\big)\Big]  \xi_k(x)dx\\[4mm]
&&+O(\eps _k^{4+p}|x_{2,k}-x_{1,k}|) \,\ \mbox{as} \,\ k\to\infty.
\end{array}\]
Since $(x+y_0)\cdot \nabla h(x+y_0)=p h(x+y_0)$, by  Proposition \ref{2:prop} and Step 1, we thus obtain from (\ref{1:H}) and above that
\[\arraycolsep=1.5pt\begin{array}{lll}
0=2\displaystyle\inte h(x+y_0)w\xi_0&=&2b_0\displaystyle\inte h(x+y_0)w\big(w+x\cdot \nabla w\big)+\sum ^2_{i=1}b_i\displaystyle\inte h(x+y_0)\frac{\partial w^2}{\partial x_i}\\[4mm]
&=&2b_0\Big[\displaystyle\inte h(x+y_0)w^2+\frac{1}{2}\inte h(x+y_0)\big(x\cdot \nabla w^2 \big)\Big]\\[4mm]
&=&2b_0\Big\{\displaystyle\inte h(x+y_0)w^2-\displaystyle\frac{1}{2}\inte w^2\big[ 2h(x+y_0)+x\cdot \nabla h(x+y_0) \big] \Big\}\\[4mm]
&=&-\displaystyle p b_0\inte h(x+y_0)w^2+b_0\displaystyle\inte   w^2\big[y_0\cdot \nabla h(x+y_0)\big]\\[4mm]
&=&-\displaystyle p b_0\inte h(x+y_0)w^2,
\end{array}\]
which therefore implies that $b_0=0$.

By the non-degeneracy assumption (\ref{1:H}), setting  $b_0=0$ into (\ref{5.2:AA}) then yields that $b_1=b_2=0$, and Step 2 is therefore proved.

%
%

\vskip 0.1truein

\noindent{\em  Step 3.} $\xi_0\equiv 0$ cannot occur.

Finally, let $y_k$ be a point satisfying $|\xi_k(y_k)|=\|\xi_k\|_{L^\infty(\R^2)}=1$. By the same argument as employed in proving Lemma \ref{lem2.3} in next section, applying the maximum principle to (\ref{3d:2}) yields that $|y_k|\le C$ uniformly in $k$. Therefore, we conclude that $\xi_k\to \xi_0\not\equiv 0$ uniformly on $\R^2$, which however contradicts to the fact that $\xi_0 \equiv 0$ on $\R^2$. This completes the proof of Theorem \ref{1:thmA}.
\qed

\section{Refined Spike Profiles}

In the following two sections, we shall derive the refined spike profiles of positive minimizers $u_k=u_{a_k}$ for $e(a_k)$ as $a_k \nearrow a^*$. The purpose of this section is to prove Theorem \ref{1:thmB}. Recall first that $u_k$ satisfies the
Euler-Lagrange equation (\ref{2:cond}). Under the assumptions of Proposition \ref{2:prop}, for convenience, we denote
\begin{equation}\label{2:beta}
 \eps _k=(a^*-a_k)^{\frac{1}{2+p}}>0,\,\ \alp _k:=\eps ^{2+p} _k>0 \ \ \mbox{and} \ \  \beta _k:=1+\frac{\mu _k\eps _k^2}{\lam ^2},
\end{equation}
where $\mu _k\in\R$ is the   Lagrange multiplier of the equation (\ref{2:cond}), so that
\[\alp _k \to 0 \ \ \mbox{and} \ \ \beta _k\to 0   \,\ \mbox{as} \,\ k\to\infty, \]
where (\ref{2:conmu}) is used.
In order to discuss the refined spike profiles of $u_k$  as $k\to\infty$, the key is thus to obtain the refined estimate of $\mu_k$ (equivalently $\beta _k$) in terms of $\eps_k$.

We next define
\begin{equation}\label{2:cong}
w_k(x):=\bar u_k(x)-w(x):=\frac{\sqrt{a^*}\eps _k}{\lam}u_k\Big(\frac{\eps _k}{\lam}x+x_k\Big)- w(x),\,\
\end{equation}
where $x_k$ is the unique maximum point of $u_k$, so that $w_k(x)\to 0$ uniformly in $\R^2$ by Proposition \ref{2:prop}. By applying (\ref{2:cond}), direct calculations then give that
$\bar u_k$ satisfies
\[
-\Delta \bar u_k(x)+\frac{\eps _k^2}{\lam ^2}V\big(\frac{\eps _k}{\lam}x+x_k\big)\bar u_k(x)=\frac{\mu _k\eps _k^2}{\lam ^2}\bar u_k(x)+\frac{a_k}{a^*}\bar u^3_k(x)\,\ \mbox{in}\,\ \R^2.
\]
Relating to the operator $\mathcal{L}:=-\Delta +(1-3w^2)$ in $\R^2$, we also denote the linearized operator
\[\mathcal{L}_k:=-\Delta +\big[1-\big(\bar u_k^2+\bar u_kw+w^2\big)\big]\ \  \mbox{in}\ \, \R^2,
\]
so that   $w_k$ satisfies
\begin{equation}\label{2:coneq}\arraycolsep=1.5pt\begin{array}{lll}
\mathcal{L}_kw_k(x)&=&-\alp _k\displaystyle\Big[\frac{1}{a^*}\bar u_k^3(x)+\frac{1}{\lam ^{2+p}} g\big(\frac{\eps_kx}{\lam}+x_k\big)h\big(x+\frac{\lam x_k}{\eps _k}\big)\bar u_k(x)\Big]\\[3mm]
&&+\beta _k\bar u_k(x)\,\ \mbox{in}\,\ \R^2,\quad \nabla w_k(0)=0,
\end{array}\end{equation}
where $V(x)=g(x)h(x)$ satisfies the assumptions of Proposition \ref{2:prop} and the coefficients $\alp_k>0$ and $\beta _k>0$ are as in (\ref{2:beta}). Define
\begin{equation}\label{2N:1}\arraycolsep=1.5pt\begin{array}{lll}
\mathcal{L}_k\psi_{1,k}(x)&=&-\alp _k\displaystyle\Big[\displaystyle\frac{1}{\lam ^{2+p}} g\big(\frac{\eps_kx}{\lam}+x_k\big)h\big(x+\frac{\lam x_k}{\eps _k}\big)\bar u_k(x)\\[3mm]
&&\qquad \quad \   +\displaystyle\frac{1}{a^*}\bar u_k^3(x)\Big]
 \,\ \mbox{in}\,\ \R^2, \quad \nabla \psi_{1,k}(0)=0,\\[4mm]
\mathcal{L}_k\psi_{2,k}(x)&=&\beta _k\bar u_k(x) \,\ \mbox{in}\,\ \R^2, \quad \nabla \psi_{2,k}(0)=0,
\end{array}\end{equation}
so that the solution $w_k(x)$ of (\ref{2:coneq}) satisfies
\begin{equation}\label{2N:3}
w_k(x):=\psi_{1,k}(x)+\psi_{2,k}(x) \,\ \mbox{in}\,\ \R^2.
\end{equation}
We first employ Proposition \ref{2:prop} to address the following estimates of $w_k$ as $k\to\infty$.

\begin{lem}\label{lem2.3} Under the assumptions of Proposition \ref{2:prop}, where $V(x)=g(x)h(x)$, we have
\begin{enumerate}
\item  $\psi_{1,k}(x)$ satisfies
  \begin{equation}\label{lem2.1:1}
  \psi_{1,k}(x)= \alp _k\psi_{1}(x)+o(\alp _k)\,\ \mbox{as} \,\ k\to\infty,
   \end{equation}
   where $\psi_{1}(x)\in C^2(\R^2)\cap L^\infty(\R^2)$ solves uniquely
     \begin{equation}\label{lem2.1:2}
 \nabla \psi_{1}(0)=0,\quad \mathcal{L}\psi_{1}(x)=- \frac{1}{a^*}w^3(x)-\frac{g(0)}{\lam ^{2+p}} h(x+y_0)w(x)
 \,\ \mbox{in}\,\ \R^2,
   \end{equation}
where $y_0\in\R^2$ is given by (\ref{1:HH}).

\item
$\psi_{2,k}(x)$ satisfies
  \begin{equation}\label{lem2.1:3}
  \psi_{2,k}(x)= \beta _k \psi_{2}(x)+o(\beta _k)\,\ \mbox{as} \,\ k\to\infty,
   \end{equation}
   where  $\psi_{2}(x)$ solves uniquely
     \begin{equation}\label{lem2.1:4}
 \nabla \psi_{2}(0)=0,\quad \mathcal{L}\psi_{2}(x)=w(x)
 \,\ \mbox{in}\,\ \R^2,
   \end{equation}
   $i.e.,$ $\psi_{2}(x)\in C^2(\R^2)\cap L^\infty(\R^2)$ satisfies
    \begin{equation}\label{lem2.1:4a}
   \psi_2=-\frac{1}{2}\big(w+x\cdot \nabla w\big).
   \end{equation}

\item $w_k$ satisfies
   \begin{equation}\label{2N:a3}
w_k(x):=\alp _k\psi_{1}(x)+\beta _k\psi_{2}(x)+o(\alp _k+\beta _k) \,\ \mbox{as} \,\ k\to\infty.
\end{equation}
\end{enumerate}
\end{lem}

\noindent\textbf{Proof.} 1.\,  We first derive $|\psi_{1,k}|\le C\alp _k$ in $\R^2$ by contradiction. On the contrary, we assume that
\begin{equation}\label{2:b2}
 \lim_{k\to\infty}\frac{\|\psi_{1,k}\|_{L^\infty}}{\alp _k}=\infty.
\end{equation}
Set $\bar \psi_{1,k}=\frac{\psi_{1,k}}{\|\psi_{1,k}\|_{L^\infty}}$ so that $\|\bar \psi_{1,k}\|_{L^\infty}=1$. Following (\ref{2N:1}), $\bar \psi_{1,k}$ then satisfies
\begin{equation}\label{2:b3}\arraycolsep=1.5pt\begin{array}{lll}
&& -\Delta \bar \psi_{1,k}  +\displaystyle \big[1-\big(\bar u_k^2+\bar u_kw+w^2\big)\big]
\bar \psi_{1,k}  \\[2mm]
&=&-\displaystyle\frac{\alp _k}{\|\psi_{1,k}\|_\infty} \Big[\displaystyle\frac{1}{\lam ^{2+p}} g\big(\frac{\eps_kx}{\lam}+x_k\big)h\big(x+\frac{\lam x_k}{\eps _k}\big)\bar u_k(x)+\displaystyle\frac{1}{a^*}\bar u_k^3(x)\Big] \,\ \mbox{in}\,\ \R^2.
\end{array}\end{equation}
Let $y_k$ be the global maximum point of $\bar \psi_{1,k}$ so that $\bar \psi_{1,k}(y_k)=\max_{y\in \R^2}\frac{\psi_{1,k}(y)}{\|\psi_{1,k}\|_{L^\infty}}=1$. Since both $\bar u_k$ and $w$ decay exponentially in view of (\ref{2:conexp}), using the maximum principle we derive from (\ref{2:b3}) that $|y_k|\le C$ uniformly in $k$.

On the other hand, applying the usual elliptic regularity theory, there exists a subsequence, still denoted by $\{\bar \psi_{1,k}\}$, of  $\{\bar \psi_{1,k}\}$ such that $\bar \psi_{1,k} \to \bar \psi _1$ weakly in $H^1(\R^2)$ and strongly in $L^q_{loc}(\R^2)$ for all $q\in [2,\infty)$. Here $\bar \psi _1$ satisfies
\[
\nabla \bar \psi_{1}(0)=0,\quad \mathcal{L}\bar \psi_{1}(x)=0  \,\ \mbox{in}\,\ \R^2,
\]
which implies that $\bar \psi_{1}=\sum ^2_{i=1}c_i\frac{\partial w}{\partial y_i}$.
Since $\nabla \bar \psi_{1}(0)=0$, we obtain that $c_1=c_2=0$. Thus, we have $\bar \psi_{1}(y)\equiv 0$ in $\R^2$, which however contradicts to the fact that $1=\bar \psi_{1,k}(y_k)\to \bar \psi_{1}(\bar y_0)$ for some $\bar y_0\in \R^2$ by passing to a subsequence if necessary. Therefore, we have $|\psi_{1,k}|\le C\alp _k$ in $\R^2$.

We next set $\phi _{1,k}(x)=\psi_{1,k}(x)-\alp _k\psi_{1}(x)$, where $\psi_{1}(x)\in C^2(\R^2)\cap L^\infty(\R^2)$ is a solution of (\ref{lem2.1:2}). Then either $\phi _{1,k}(x)=O(\alp _k)$ or $\phi _{1,k}(x)=o(\alp _k)$ as $k\to\infty$, and $\phi _{1,k}$ satisfies
\[
\nabla \phi _{1,k}(0)=0,\quad -\Delta \phi_{1,k}  + \big[1-\big(\bar u_k^2+\bar u_kw+w^2\big)\big]\phi_{1,k}=-\alp _kf_k(x)  \,\ \mbox{in}\,\ \R^2,
\]
where $f_k(x)$ satisfies
\[\arraycolsep=1.5pt\begin{array}{lll}
f_k(x)&=&\big(2w^2-\bar u_k^2-\bar u_kw\big)\psi_{1}(x)+\displaystyle\frac{1}{a^*}\big(\bar u_k^3(x)-w^3(x)\big)\\[3mm]
&&+\displaystyle\frac{1}{\lam ^{2+p}}\Big[g\big(\frac{\eps_kx}{\lam}+x_k\big)h\big(x+\frac{\lam x_k}{\eps _k}\big)\bar u_k(x)-g(0)h(x+y_0)w(x)\Big].
\end{array}\]
One can note that $f_k(x)\to 0$ uniformly as $k\to\infty$. Therefore, applying the previous argument yields necessarily that $\phi _{1,k}(x)=o(\alp _k)$ as $k\to\infty$, and the proof of (\ref{lem2.1:1}) is then complete. Also, the property (\ref{2:linearized}) gives the uniqueness of solutions for (\ref{lem2.1:2}).

2.\,   Since the proof of (\ref{lem2.1:3}) is very similar to that of (\ref{lem2.1:1}), we omit the details. Further, the property (\ref{2:linearized}) gives the uniqueness of $\psi_2$. Also, one can check directly that $-(w+x\cdot \nabla w)/2$ is a solution of (\ref{lem2.1:4}), which therefore implies that  (\ref{lem2.1:4a}) holds.

3.\, The expansion (\ref{2N:a3}) now follows immediately from (\ref{2N:3}), (\ref{lem2.1:1}) and (\ref{lem2.1:3}), and the proof is therefore complete.
\qed

\subsection{Proof of Theorem \ref{1:thmB}}

The main aim of this subsection is to prove Theorem \ref{1:thmB} on the refined spike behavior of positive minimizers. In this whole subsection, we assume that the potential $V(x)=h(x)\in C^2(\R^2)$ satisfies $\lim_{|x|\to\infty} h(x) = \infty$ and (\ref{1:H}), where $h(x)$ is homogeneous of degree $p\ge 2$. Following (\ref{2:beta}), from now on we denote for simplicity that
\begin{equation}\label{3ab:20}
o\big([\alp _k+\beta _k]^2\big)=o\big(\alp _k^2\big)+o\big( \alp _k \beta _k \big)+o\big(\beta _k^2\big)\,\ \mbox{as} \,\ k\to\infty,
\end{equation}
where $\alp _k$ and $\beta_k$ are defined in (\ref{2:beta}). We first use Lemma \ref{lem2.3} to establish the following lemmas.

\begin{lem}\label{lem3.1} Suppose that $V(x)=h(x)\in C^2(\R^2)$ satisfies $\lim_{|x|\to\infty} h(x) = \infty$ and (\ref{1:H}) for some $y_0\in \R^2$, where $h(x)$ is homogeneous of degree $p\ge 2$.  Then there exists an $x_0\in\R^2$ such that the unique maximum point $x_k$ of $u_k$   satisfies
\begin{equation}\label{3a:20}
\Big|\alp _k\Big(\frac{\lam x_k}{\eps _k}-y_0\Big)-\alp_k\beta _k\frac{y_0}{2}\Big|= \alp _k^2\, O (|x_0|) +o\big([\alp _k+\beta _k]^2\big)\,\ \mbox{as} \,\ k\to\infty.
\end{equation}
\end{lem}

\noindent\textbf{Proof.}
Multiplying (\ref{lem2.1:2}) and (\ref{lem2.1:4}) by $\frac{\partial w}{\partial x_1}$ and then integrating over $\R^2$, respectively, we obtain from  (\ref{1:H}) and (\ref{2:linearized}) that
\begin{equation}\label{3a:1a}
\inte \frac{\partial w}{\partial x_1}\mathcal{L}w_k =\inte\frac{\partial w}{\partial x_1}w=\inte\frac{\partial w}{\partial x_1}h(x+y_0)w=0,
\end{equation}
where $y_0$ is given by the assumption (\ref{1:H}). Similarly, we derive from (\ref{2:coneq}) and (\ref{2N:a3}) that
\begin{equation}\label{3a:21}\arraycolsep=1.5pt\begin{array}{lll}
 \displaystyle\inte\frac{\partial w}{\partial x_1}\mathcal{L}_kw_k
&=&\displaystyle\beta _k\inte\frac{\partial w}{\partial x_1}\bar u_k -\alp _k\displaystyle\inte\frac{\partial w}{\partial x_1}\Big[\frac{1}{a^*}\bar u_k^3 +\frac{\bar u_k }{\lam ^{2+p}} h\big(x+\frac{\lam x_k}{\eps _k}\big)\Big]\\[3mm]
&=& \displaystyle\alp _k\beta _k\inte\frac{\partial w}{\partial x_1}\psi_1 +o(\alp _k\beta _k+\beta^2 _k)- I_1,
\end{array}\end{equation}
where the identity $\inte\frac{\partial w}{\partial x_1}\psi_2 =0$ is used, since $\frac{\partial w}{\partial x_1}\psi_2$ is odd in $x_1$ by the radial symmetry of $\psi_2$.
We obtain from  (\ref{1:H}) and (\ref{3a:1a}) that
\begin{equation}\label{3a:23}\arraycolsep=1.5pt\begin{array}{lll}
I_1&=& \alp _k\displaystyle\inte\frac{\partial w}{\partial x_1}\Big[\frac{1}{a^*}\bar u_k^3 +\frac{\bar u_k }{\lam ^{2+p}} h\big(x+\frac{\lam x_k}{\eps _k}\big)\Big] \\[4mm]
&=& \alp _k\displaystyle\inte\frac{\partial w}{\partial x_1}\Big\{\frac{1}{a^*}\big(\bar u_k^3 -w^3\big)+\displaystyle\frac{1}{\lam ^{2+p}} \Big[h\big(x+\frac{\lam x_k}{\eps _k}\big)\bar u_k -h(x+y_0)w\Big]\Big\}\\[4mm]
&=& \displaystyle\frac{\alp _k}{a^*}\inte\frac{\partial w}{\partial x_1}w_k\big(3w^2+3ww_k+w^2_k\big)\\[4mm]
&&+\displaystyle\frac{\alp _k}{\lam ^{2+p}} \inte\frac{\partial w}{\partial x_1}\Big[h\big(x+\frac{\lam x_k}{\eps _k}\big)\bar u_k-h(x+y_0)w\Big]\\[4mm]
&=&\displaystyle\frac{3\alp_k^2}{a^*}\inte\frac{\partial w}{\partial x_1}w^2\psi_1+o(\alp _k^2+\alp_k\beta _k)+I_2,
\end{array}\end{equation}
where we have used the identity $\inte\frac{\partial w}{\partial x_1}w^2\psi_2 =0$, since $\frac{\partial w}{\partial x_1}w^2\psi_2$ is odd in $x_1$ by the radial symmetry of $\psi_2$. Further, applying (\ref{2N:a3}) and (\ref{3a:1a}) yields that
\begin{equation}\label{3a:24}\arraycolsep=1.5pt\begin{array}{lll}
\displaystyle\frac{\lam ^{2+p}}{\alp _k} I_2&=&\displaystyle\inte\frac{\partial w}{\partial x_1}\Big\{h\big(x+\frac{\lam x_k}{\eps _k}\big)\big[\bar u_k-w\big]+\big[h\big(x+\frac{\lam x_k}{\eps _k}\big)-h(x+y_0)\big]w\Big\}\\[4mm]
&=&\displaystyle\inte\frac{\partial w}{\partial x_1}h(x+y_0)w_k+o(\alp _k+\beta _k)\\[4mm]
&&+\displaystyle\inte\frac{\partial w}{\partial x_1} \Big[\Big(\frac{\lam x_k}{\eps _k}-y_0\Big)\cdot \nabla h(x+y_0)\Big]w +o\Big(\Big|\frac{\lam x_k}{\eps _k}-y_0\Big|\Big)\\[4mm]
&=&\alp _k\displaystyle\inte\frac{\partial w}{\partial x_1}h(x+y_0)\psi_1+\beta_k\displaystyle\inte\frac{\partial w}{\partial x_1}h(x+y_0)\psi_2\\[4mm]
&&+\displaystyle\inte\frac{\partial w}{\partial x_1}\Big[\Big(\frac{\lam x_k}{\eps _k}-y_0\Big)\cdot \nabla h(x+y_0)\Big]w+o\Big(\alp _k+\big|\frac{\lam x_k}{\eps _k}-y_0\big|+\beta _k \Big),
\end{array}\end{equation}
where  (\ref{2:conb}) is used for the second identity. We thus get that
\begin{equation}\label{3a:25}\arraycolsep=1.5pt\begin{array}{lll}
I_1&=&\alp^2_k\Big[\displaystyle\frac{3}{a^*}\inte\frac{\partial w}{\partial x_1}w^2\psi_1+\displaystyle\frac{1}{\lam ^{2+p}}\inte\frac{\partial w}{\partial x_1}h(x+y_0)\psi _1\Big]\\[4mm]
&&+\displaystyle\frac{\alp _k\beta _k}{\lam ^{2+p}}\inte\frac{\partial w}{\partial x_1}h(x+y_0)\psi_2+\displaystyle\frac{\alp _k}{\lam ^{2+p}}\inte\frac{\partial w}{\partial x_1} \Big[\Big(\frac{\lam x_k}{\eps _k}-y_0\Big)\cdot \nabla h(x+y_0)\Big]w\\[4mm]
&&+\displaystyle o\Big(\alp_k\big|\frac{\lam x_k}{\eps _k}-y_0\big|+[\alp _k+\beta _k]^2\Big).
\end{array}\end{equation}

On the other hand, we obtain from  (\ref{3a:1a}) that
\begin{equation}\label{3a:22}\arraycolsep=1.5pt\begin{array}{lll}
\displaystyle\inte\frac{\partial w}{\partial x_1}\mathcal{L}_kw_k&=& \displaystyle\inte\frac{\partial w}{\partial x_1}\mathcal{L}w_k +\displaystyle\inte\frac{\partial w}{\partial x_1}\big(\mathcal{L}_k-\mathcal{L}\big)w_k \\[3mm]
&=&  -\displaystyle\inte\frac{\partial w}{\partial x_1}w_k^2(3w+w_k)\\[3mm]
&=&-3\alp _k^2\displaystyle\inte\frac{\partial w}{\partial x_1}w\psi _1^2-6\alp_k\beta_k\displaystyle\inte\frac{\partial w}{\partial x_1}w\psi _1\psi _2+o(\alp _k^2+\alp_k\beta _k).
\end{array}\end{equation}
Combining (\ref{3a:21}), (\ref{3a:22}) and (\ref{3a:25}), we now conclude from (\ref{1:H}) and (\ref{2N:a3})  that
\begin{equation}\label{3a:26}\arraycolsep=1.5pt\begin{array}{lll}
&& \displaystyle\frac{\alp _k}{\lam ^{2+p}}\inte\frac{\partial w}{\partial x_1} \Big[\Big(\frac{\lam x_k}{\eps _k}-y_0\Big)\cdot \nabla h(x+y_0)\Big]w\\[4mm]
&=&
\alp_k\beta _k\Big[\displaystyle\inte\frac{\partial w}{\partial x_1}\psi_1+6\displaystyle\inte\frac{\partial w}{\partial x_1}w\psi_1\psi_2-\displaystyle\frac{1}{\lam ^{2+p}}\inte\frac{\partial w}{\partial x_1}h(x+y_0)\psi_2\Big]\\[4mm]
&&-\alp _k^2\Big[\displaystyle\frac{3}{a^*}\inte\frac{\partial w}{\partial x_1}w^2\psi_1+\displaystyle\frac{1}{\lam ^{2+p}}\inte\frac{\partial w}{\partial x_1}h(x+y_0)\psi _1-3\displaystyle\inte\frac{\partial w}{\partial x_1}w\psi _1^2\Big]\\[3mm]
&&+o([\alp _k+\beta _k]^2).
\end{array}\end{equation}
We {\em claim} that the coefficient $I_3$ of the term $\alp_k\beta_k$ in (\ref{3a:26}) satisfies
\begin{equation}\label{3a:26a}\arraycolsep=1.5pt\begin{array}{lll}
I_3:&=&\displaystyle\inte\frac{\partial w}{\partial x_1}\psi_1+6\displaystyle\inte\frac{\partial w}{\partial x_1}w\psi_1\psi_2-\displaystyle\frac{1}{\lam ^{2+p}}\inte\frac{\partial w}{\partial x_1}h(x+y_0)\psi_2 \\[3mm]
&=&\displaystyle\frac{1}{2\lam ^{2+p}}\inte w\Big[y_0\cdot \nabla h(x+y_0)\Big]\frac{\partial w}{\partial x_1}.
\end{array}\end{equation}
If  (\ref{3a:26a}) holds, we then derive from (\ref{3a:26}) that there exists some $x_0=(x_{10},x_{20})\in\R^2$ such that
\begin{equation}\label{3M:26}\arraycolsep=1.5pt\begin{array}{lll}
&& \displaystyle\frac{1}{2\lam ^{2+p}}\inte\frac{\partial w^2}{\partial x_j} \Big[\alp _k\Big(\frac{\lam x_k}{\eps _k}-y_0\Big)-\alp_k\beta _k\frac{y_0}{2}\Big]\cdot \nabla h(x+y_0) \\[4mm]
&=&\alp _k^2O(|x_{j0}|)+o([\alp _k+\beta _k]^2),\quad j=1,2.
\end{array}\end{equation}
By the non-degeneracy assumption of (\ref{1:H}), we further conclude from (\ref{3M:26})  that (\ref{3a:20}) holds for some $x_0\in\R^2$, and the lemma is therefore proved.

To complete the proof of the lemma, the rest is to prove the claim  (\ref{3a:26a}).
Indeed, using the integration by parts, we derive from (\ref{lem2.1:4a}) that
\[\arraycolsep=1.5pt\begin{array}{lll}
 A:&=&\displaystyle\inte\frac{\partial w}{\partial x_1}\psi_1+6\displaystyle\inte\frac{\partial w}{\partial x_1}w\psi_1\psi_2\\[4mm]
 &=&\displaystyle\inte\frac{\partial w}{\partial x_1}\psi_1-3\displaystyle\inte\frac{\partial w}{\partial x_1}w^2\psi_1- \displaystyle\frac 3 2 \inte\frac{\partial w}{\partial x_1}\psi_1 (x\cdot \nabla w^2)\\[4mm]
 &=&\displaystyle\inte\frac{\partial w}{\partial x_1}\psi_1-3\displaystyle\inte\frac{\partial w}{\partial x_1}w^2\psi_1
 \\[4mm]&&
 + \displaystyle\frac 3 2 \inte w^2\Big[2\frac{\partial w}{\partial x_1}\psi_1+x\cdot \nabla \Big(\frac{\partial w}{\partial x_1}\psi_1\Big)\Big]\\[4mm]
 &=&\displaystyle\inte\frac{\partial w}{\partial x_1}\psi_1+\displaystyle\frac 3 2 \inte w^2\Big[\frac{\partial w}{\partial x_1}(x\cdot \nabla \psi_1)+\psi_1 x\cdot \nabla \Big(\frac{\partial w}{\partial x_1}\Big)\Big].
\end{array}\]
Since $(x+y_0)\cdot \nabla h(x+y_0)=ph(x+y_0)$, we obtain from (\ref{1:H}), (\ref{lem2.1:4a}) and (\ref{3a:1a}) that
\[\arraycolsep=1.5pt\begin{array}{lll}
 B:&=&-\displaystyle\frac{1}{\lam ^{2+p}}\inte\frac{\partial w}{\partial x_1}h(x+y_0)\psi_2\\[4mm]
 &=&\displaystyle\frac{1}{2\lam ^{2+p}}\inte\frac{\partial w}{\partial x_1}h(x+y_0)(w+x\cdot \nabla w)\\[4mm]
 &=&\ - \displaystyle\frac{1}{2\lam ^{2+p}}\inte w\Big[2h(x+y_0)\frac{\partial w}{\partial x_1}+x\cdot \nabla \Big(\frac{\partial w}{\partial x_1}h(x+y_0)\Big)\Big]\\[4mm]
 &=& - \displaystyle\frac{1}{2\lam ^{2+p}}\inte w\Big\{\big[x\cdot \nabla h(x+y_0)\big]\frac{\partial w}{\partial x_1}+h(x+y_0)x\cdot \nabla \Big(\frac{\partial w}{\partial x_1}\Big)\Big\}\\[4mm]
 &=& - \displaystyle\frac{1}{2\lam ^{2+p}}\inte wh(x+y_0)\Big[x\cdot \nabla \Big(\frac{\partial w}{\partial x_1}\Big)\Big]
 +\displaystyle\frac{1}{2\lam ^{2+p}}\inte w\Big[y_0\cdot \nabla h(x+y_0)\Big]\frac{\partial w}{\partial x_1}.
\end{array}\]
By above calculations, we then get from (\ref{3a:26a}) that
\begin{equation}\label{3a:26b}\arraycolsep=1.5pt\begin{array}{lll}
I_3=A+B&=&\displaystyle\inte\frac{\partial w}{\partial x_1}\psi_1+\displaystyle\frac 3 2 \inte w^2 \frac{\partial w}{\partial x_1}(x\cdot \nabla \psi_1)\\[4mm]
&&+\displaystyle\frac 1 2 \inte \Big[3w^2\psi_1-\frac{wh(x+y_0)}{ \lam ^{2+p}}\Big]\Big[x\cdot \nabla \Big(\frac{\partial w}{\partial x_1}\Big)\Big]
\\[4mm]
&&
 +\displaystyle\frac{1}{2\lam ^{2+p}}\inte w\Big[y_0\cdot \nabla h(x+y_0)\Big]\frac{\partial w}{\partial x_1}\\[4mm]
&:=&
 I_4+\displaystyle\frac{1}{2\lam ^{2+p}}\inte w\Big[y_0\cdot \nabla h(x+y_0)\Big]\frac{\partial w}{\partial x_1}.
\end{array}\end{equation}

Applying the integration by parts, we derive from  (\ref{lem2.1:2}) that
\[
\arraycolsep=1.5pt\begin{array}{lll}
&&\displaystyle\inte\frac{\partial w}{\partial x_1}\psi_1+\displaystyle\frac 1 2 \inte \Big[3w^2\psi_1-\frac{wh(x+y_0)}{ \lam ^{2+p}}\Big]\Big[x\cdot \nabla \Big(\frac{\partial w}{\partial x_1}\Big)\Big]\\[4mm]
&=&\displaystyle\inte\frac{\partial w}{\partial x_1}\psi_1+\displaystyle\frac 1 2 \inte \Big[3w^2\psi_1-\Big(\frac{w^3}{a^*}+\frac{wh(x+y_0)}{ \lam ^{2+p}}\Big)\Big]\Big[x\cdot \nabla \Big(\frac{\partial w}{\partial x_1}\Big)\Big]\\[4mm]
&=&\displaystyle\inte\frac{\partial w}{\partial x_1}\psi_1+\displaystyle\frac 1 2 \inte \big(-\Delta \psi_1+\psi_1\big)\Big[x\cdot \nabla \Big(\frac{\partial w}{\partial x_1}\Big)\Big]\\[4mm]
&=&\displaystyle\inte\frac{\partial w}{\partial x_1}\psi_1+\displaystyle\frac 1 2 \inte \big(-\Delta \psi_1\big)\Big[x\cdot \nabla \Big(\frac{\partial w}{\partial x_1}\Big)\Big]-\displaystyle\frac 1 2 \inte\frac{\partial w}{\partial x_1}\big[2\psi_1+x\cdot \nabla \psi_1\big]\\[4mm]
&=&\displaystyle\frac 1 2 \inte \big(-\Delta \psi_1\big)\Big[x\cdot \nabla \Big(\frac{\partial w}{\partial x_1}\Big)\Big]-\displaystyle\frac 1 2 \inte\frac{\partial w}{\partial x_1}\big(x\cdot \nabla \psi_1\big),
\end{array}\]
which then gives from (\ref{3a:26b}) that
\begin{equation}\label{3a:26c}\arraycolsep=1.5pt\begin{array}{lll}
-2I_4&=& \displaystyle  \inte  \Delta \psi_1 \Big[x\cdot \nabla \Big(\frac{\partial w}{\partial x_1}\Big)\Big]+\displaystyle \inte\frac{\partial }{\partial x_1}\Big[w-w^3\Big]\big(x\cdot \nabla \psi_1\big) \\[4mm]
&=&\displaystyle  \inte  \Delta \psi_1 \Big[x\cdot \nabla \Big(\frac{\partial w}{\partial x_1}\Big)\Big]+\displaystyle  \inte\frac{\partial \Delta w}{\partial x_1}  \big(x\cdot \nabla \psi_1\big).
\end{array}\end{equation}
To further simplify $I_4$, we next rewrite $\psi_1$ as $\psi_1(x)=\psi_1(r,\theta )$, where $x=r(\cos \theta ,\sin \theta)$ and $(r,\theta )$ is the polar coordinate in $\R^2$. We then follow from (\ref{lem2.1:2}) and (\ref{3a:26c}) that
\begin{equation}\label{3a:26d}\arraycolsep=1.5pt\begin{array}{lll}
-2I_4&=& \displaystyle  \intr\int ^{2\pi}_0 \Big\{\big[r\big(\psi_1\big)_r\big]_r+\frac{1}{r}\big(\psi_1\big)_{\theta\theta}\Big\}r\frac{\partial }{\partial r}
\big(w'\cos \theta\big)d\theta dr\\[4mm]
&&+\displaystyle  \intr\int ^{2\pi}_0 \frac{\partial }{\partial r}\Big(w''+\frac{w'}{r}\Big)\cos\theta \, r^2\big(\psi_1\big)_rd\theta dr\\[4mm]
&=&-\displaystyle  \intr\int ^{2\pi}_0 r\big(\psi_1\big)_r(rw'')'\cos \theta d\theta dr\\[4mm]
&&+\displaystyle  \intr\int ^{2\pi}_0 r\big(\psi_1\big)_r\Big[r\frac{\partial }{\partial r}\Big(w''+\frac{w'}{r}\Big)\Big]\cos\theta   d\theta dr\\[4mm]
&&+\displaystyle  \intr\int ^{2\pi}_0\big(\psi_1\big)_{\theta\theta}w''\cos \theta d\theta dr\\[4mm]
&=&-\displaystyle  \intr\int ^{2\pi}_0 r\big(\psi_1\big)_r\Big\{(rw'')'-\Big[r\frac{\partial }{\partial r}\Big(w''+\frac{w'}{r}\Big)\Big]\Big\}\cos \theta d\theta dr\\[4mm]
&&-\displaystyle  \intr\int ^{2\pi}_0 \psi_1 w''\cos \theta d\theta dr\\[4mm]
&=&-\displaystyle  \intr\int ^{2\pi}_0 \big(\psi_1\big)_rw'\cos \theta d\theta dr-\displaystyle  \intr\int ^{2\pi}_0 \psi_1 w''\cos \theta d\theta dr=0,
\end{array}\end{equation}
$i.e.,$ $I_4=0$, which  therefore implies that the claim  (\ref{3a:26a}) holds by applying (\ref{3a:26b}).
\qed

\begin{rem}\label{rem3.1}  Whether the point $x_0\in\R^2$ in Lemma \ref{lem3.1} is the origin or not is determined completely by the fact that whether the coefficient $I_5$ of the term $\alp_k^2$ in (\ref{3a:26}) is zero or not, where $I_5$ satisfies
\[
I_5:=\displaystyle\frac{3}{a^*}\inte\frac{\partial w}{\partial x_1}w^2\psi_1+\displaystyle\frac{1}{\lam ^{2+p}}\inte\frac{\partial w}{\partial x_1}h(x+y_0)\psi _1-3\displaystyle\inte\frac{\partial w}{\partial x_1}w\psi _1^2.
\]
If $h(x)$ is not even in $x$, it however seems difficult to derive that whether $I_5  =0$ or not.
\end{rem}

\begin{lem}\label{lem3.2} Suppose that $V(x)=h(x)\in C^2(\R^2)$ satisfies $\lim_{|x|\to\infty} h(x) = \infty$ and (\ref{1:H}) for some $y_0\in\R^2$, where $h(x)$ is homogeneous of degree $p\ge 2$. Then we have
\begin{equation}\label{lem3.2:1}
w_k :=\alp _k\psi_{1} +\beta _k\psi_{2} + \alp ^2_k\psi_{3} +\beta ^2_k\psi_{4} +\alp _k\beta _k\psi_{5}+o([\alp_k+\beta _k]^2) \ \ \mbox{as} \ \ k\to\infty,
\end{equation}
where $\psi_1(x),\, \psi_2(x)  \in C^2(\R^2)\cap L^\infty(\R^2)$ are given in Lemma \ref{lem2.3} with $g(0)=1$, and $\psi_{i}(x)  \in C^2(\R^2)\cap L^\infty(\R^2)$, $i=3,4,5$, solves uniquely
     \begin{equation}\label{lem3.2:2}
 \nabla \psi_{i}(0)=0\ \, \mbox{and}\,\ \mathcal{L}\psi_{i}(x)=f_i(x)
 \,\ \mbox{in}\,\ \R^2,\,\ i=3,4,5,
   \end{equation}
 and $f_i(x)$ satisfies for some $y^0\in\R^2$,
 \begin{equation}\label{lem3.2:3}\arraycolsep=1.5pt
 f_i(x)=\left\{\begin{array}{lll}
 &3w\psi^2_1-\displaystyle\Big(\frac{3w^2}{a^*}+\frac{h(x+y_0)}{\lam ^{2+p}}\Big)\psi_1-\displaystyle\frac{w}{\lam ^{1+p}}\big[y^0\cdot \nabla h(x+y_0)\big],\,\ &\mbox{if}\ \ i=3;\\[3mm]
 &3w\psi^2_2+\psi_2, \,\ &\mbox{if}\ \ i=4;\\
 &  6w\psi_1\psi_2+\psi _1-\displaystyle\Big(\frac{3w^2}{a^*}+\displaystyle\frac{h(x+y_0)}{\lam ^{2+p}}\Big)\psi_2 \\[2mm]
 &-\displaystyle\frac{w}{2\lam ^{2+p}} \Big[y_0\cdot \nabla h(x+y_0)\Big] , \,\ &\mbox{if}\ \ i=5;
 \end{array}\right.
    \end{equation}
where $y_0\in\R^2$ is given by (\ref{1:H}). Moreover, $\psi_{4}$  is radially symmetric.
\end{lem}

\noindent\textbf{Proof.} Following Lemma \ref{lem2.3}(3), set
\[
v_k=w_k-\alp_k\psi_1-\beta_k\psi_2,
\]
so that
\begin{equation}\label{3a:31}\arraycolsep=1.5pt\begin{array}{lll}
 \mathcal{L}_kw_k&=& \mathcal{L}_k(v_k+\alp_k\psi_1+\beta_k\psi_2) \\[1mm]
&=& \mathcal{L}_kv_k+\alp_k(\mathcal{L}_k-\mathcal{L})\psi_1+\beta_k(\mathcal{L}_k-
\mathcal{L})\psi_2+\alp_k\mathcal{L}\psi_1+\beta_k\mathcal{L}\psi_2\\ [1mm]
&=&\mathcal{L}_kv_k-w_k(\alp_k\psi_1+\beta_k\psi_2)(3w+w_k)-
\alp_k\Big[\displaystyle\frac{w^3}{a^*}+\displaystyle\frac{h(x+y_0)w}{\lam ^{2+p}}\Big]+\beta_kw.
\end{array}\end{equation}
Applying (\ref{2:coneq}), we then have
\begin{equation}\label{3a:31a}\arraycolsep=1.5pt\begin{array}{lll}
 \mathcal{L}_kv_k&=&\mathcal{L}_kw_k+w_k(\alp_k\psi_1+\beta_k\psi_2)(3w+w_k)+
\alp_k\Big[\displaystyle\frac{w^3}{a^*}+\displaystyle\frac{h(x+y_0)w}{\lam ^{2+p}}\Big]-\beta_kw\\[2mm]
&=&w_k(\alp_k\psi_1+\beta_k\psi_2)(3w+w_k)+\beta_k(\bar u_k-w)\\[1mm]
&&- \alp_k\Big\{\displaystyle\frac{1}{a^*}(\bar u_k^3-w^3)+\displaystyle\frac{1}{\lam ^{2+p}}\Big[h\big(x+\frac{\lam x_k}{\eps _k}\big)\bar u_k-h(x+y_0)w \Big]\Big\}\\[2mm]
&=&w_k(\alp_k\psi_1+\beta_k\psi_2)(3w+w_k)+\beta_kw_k-I_6,
\end{array}\end{equation}
where $I_6$ satisfies
\[\arraycolsep=1.5pt\begin{array}{lll}
 I_6&=&\displaystyle\frac{\alp_k}{a^*}w_k(3w^2+3ww_k+w_k^2)\\[1mm]
 &&+\displaystyle\frac{\alp_k}{\lam ^{2+p}}\Big\{h(x+y_0)(\bar u_k-w)+\Big[h\big(x+\frac{\lam x_k}{\eps _k}\big)-h(x+y_0)\Big]\bar u_k \Big\}\\[2mm]
&=& \displaystyle\frac{\alp_k}{a^*}w_k(3w^2+3ww_k+w_k^2)+\displaystyle\frac{\alp_k}{\lam ^{2+p}}h(x+y_0)w_k
\\[2mm]
&&+\displaystyle\frac{\alp_k}{\lam ^{2+p}}\Big[\Big(\frac{\lam x_k}{\eps _k}-y_0\Big)\cdot \nabla h(x+y_0)\Big]\bar u_k +o\big([\alp _k+\beta _k]^2\big)\\[3mm]
&=&\alp_kw_k\displaystyle\Big(\frac{3w^2}{a^*}+\frac{h(x+y_0)}{\lam ^{2+p}}\Big)+\displaystyle\frac{\alp_k}{a^*}w_k^2(3w+w_k)\\[3mm]
&&+\displaystyle\frac{\alp_k}{\lam ^{2+p}}\Big[\Big(\frac{\lam x_k}{\eps _k}-y_0\Big)\cdot \nabla h(x+y_0)\Big]\bar u_k +o\big([\alp _k+\beta _k]^2\big),
\end{array}\]
where Lemma \ref{lem3.1} is used in the second equality.
By Lemma \ref{lem3.1} again, there exists $y^0\in\R^2$ such that
\[
 \Big|\alp _k\Big(\frac{\lam x_k}{\eps _k}-y_0\Big)-\alp_k\beta _k\frac{y_0}{2}-\alp _k^2y^0\Big|=o\big([\alp _k+\beta _k]^2\big) \,\ \mbox{as}\,\ k\to\infty.
\]
We thus obtain from above that
\begin{equation}\label{3a:33}\arraycolsep=1.5pt\begin{array}{lll}
 \mathcal{L}_kv_k&=&w_k(\alp_k\psi_1+\beta_k\psi_2)(3w+w_k)+\beta_kw_k-
 \alp_kw_k\displaystyle\Big(\frac{3w^2}{a^*}+\frac{h(x+y_0)}{\lam ^{2+p}}\Big)\\[2mm]
 &&-\displaystyle\frac{\alp_k}{\lam ^{2+p}}\Big[\Big(\frac{\lam x_k}{\eps _k}-y_0\Big)\cdot \nabla h(x+y_0)\Big]\bar u_k-\displaystyle\frac{\alp_k}{a^*}w_k^2(3w+w_k)+o([\alp_k+\beta _k]^2)\\[2mm]
 &=&\alp_k^2\displaystyle\Big\{3w\psi^2_1-\Big(\frac{3w^2}{a^*}+\frac{h(x+y_0)}{\lam ^{2+p}}\Big)\psi_1 -\displaystyle\frac{w}{\lam ^{1+p}}\Big[y^0\cdot \nabla h(x+y_0)\Big]\Big\}\\[2mm]
 &&+\alp_k\beta_k\displaystyle\Big\{ 6w\psi_1\psi_2+\psi _1-\Big(\frac{3w^2}{a^*}+\frac{h(x+y_0)}{\lam ^{2+p}}\Big)\psi_2 -\displaystyle\frac{1}{2\lam ^{2+p}}w\Big[y_0\cdot \nabla h(x+y_0)\Big]\Big\}\\[3mm]
 &&+\beta_k^2(3w\psi^2_2+\psi_2)+o([\alp_k+\beta _k]^2)\ \, \ \mbox{in}\ \ \R^2.
\end{array}
\end{equation}
Following (\ref{3a:33}), the same argument of proving Lemma \ref{lem2.3} then gives (\ref{lem3.2:1}). Finally, since $f_4(x)$ is radially symmetric, there exists a radial solution $\psi_4$. Further, the property (\ref{2:linearized}) gives the uniqueness of $\psi_4$. Therefore, $\psi_4$ must be radially symmetric, and the proof is complete.
\qed

\begin{lem}\label{lem3.3} Suppose that $V(x)=h(x)\in C^2(\R^2)$ satisfies $\lim_{|x|\to\infty} h(x) = \infty$ and (\ref{1:H}) for some $y_0\in\R^2$, where $h(x)$ is homogeneous of degree $p\ge 2$.  Then we have
\begin{equation}\label{3a:37a}
 \displaystyle\inte w\psi_{1} =0,\ \ \inte w\psi_{2} =0,
\end{equation} and
\begin{equation}\label{3a:37aa}
I= \inte \big(2w\psi_4+\psi^2_2\big)=0.
\end{equation}
However, we have
 \begin{equation}\label{a:claim3}
II=2\inte w\psi_5+2\inte \psi_1\psi_2=-\frac{2+p}{2}<0.
\end{equation}
Here $\psi_1(x), \cdots ,\psi_5(x)\in C^2(\R^2)\cap L^\infty(\R^2)$  are given in Lemma \ref{lem2.3} with $g(0)=1$ and Lemma \ref{lem3.2}.
\end{lem}

Since the proof of Lemma \ref{lem3.3} is very involved, we leave it to the appendix. Following above lemmas, we are now ready to derive the comparison relation between $\beta_k$ and $\alp _k$.

\begin{prop}\label{prop3.4} Suppose that $V(x)=h(x)\in C^2(\R^2)$ satisfies $\lim_{|x|\to\infty} h(x) = \infty$ and (\ref{1:H}) for  some $y_0\in\R^2$, where $h(x)$ is homogeneous of degree $p\ge 2$.  Then we have
\begin{equation}\label{3:35}
\beta _k=C^*\alp_k\,\ \mbox{as} \,\ k\to\infty,
\end{equation}
where the constant $C^*$ satisfies
\begin{equation}\label{3:35a}
C^*= \frac{2}{2+p}\Big(2\displaystyle\inte w\psi_{3}+\displaystyle\inte \psi_{1}^2\Big)\not =0.
\end{equation}
Moreover, $w _k$ satisfies
\begin{equation}\label{3:36}
w_k :=\big[\psi_{1}+ C^*\psi_{2} \big] \alp_k+\big[\psi_{3}+(C^*)^2\psi_{4}+ C^* \psi_{5}\big]\alp_k^2+o(\alp_k^2) \,\ \mbox{as} \,\ k\to\infty,
\end{equation}
Here $\psi_1(x),\cdots ,\psi_5(x)\in C^2(\R^2)\cap L^\infty(\R^2)$ are given in Lemma \ref{lem2.3} with $g(0)=1$ and Lemma \ref{lem3.2}.
\end{prop}

\noindent\textbf{Proof.} Note from (\ref{2:cong}) that $w_k$ satisfies
\begin{equation}\label{2a:k}
 \inte w^2=\inte \bar u_k^2=\inte \big(w+w_k\big)^2,\  i.e. \ \ 2\inte ww_k+\inte w_k^2=0.
\end{equation}
Applying (\ref{2a:k}), we then derive from Lemma  \ref{lem3.2}   that
\begin{equation}\label{3a:37}\arraycolsep=1.5pt\begin{array}{lll}
0&=&2\displaystyle\inte ww_k+\displaystyle\inte w_k^2\\[3mm]
&=&2\displaystyle\inte w(\alp _k\psi_{1}  +\beta _k\psi_{2}+ \alp ^2_k\psi_{3} +\beta ^2_k\psi_{4} +\alp _k\beta _k\psi_{5})\\[4mm]
&&+\displaystyle\inte (\alp _k\psi_{1} +\beta _k\psi_{2} + \alp ^2_k\psi_{3} +\beta ^2_k\psi_{4} +\alp _k\beta _k\psi_{5})^2+o([\alp_k+\beta _k]^2)\\[4mm]
&=& \alp _k\Big(\displaystyle2\inte  w\psi_{1}\Big)+\beta _k\Big(\displaystyle 2\inte  w\psi_{2}\Big)+\beta _k^2\Big(2\displaystyle\inte  w\psi_{4}+\displaystyle\inte \psi_{2}^2\Big)\\[4mm]
&&+\alp_k\beta _k\Big(2\displaystyle\inte w\psi_{5}+2\displaystyle\inte \psi_{1}\psi_{2}\Big)+\alp _k^2\Big(2\displaystyle\inte w\psi_{3}+\displaystyle\inte \psi_{1}^2\Big)+o([\alp_k+\beta _k]^2)\\[4mm]
&=&-\displaystyle\frac{2+p}{2}\alp_k\beta _k+\alp _k^2\Big(2\displaystyle\inte w\psi_{3}+\displaystyle\inte \psi_{1}^2\Big)+o([\alp_k+\beta _k]^2),
\end{array}
\end{equation}
where Lemma  \ref{lem3.3} is used in the last equality. It then follows from (\ref{3a:37}) that
\[
2\displaystyle\inte w\psi_{3}+\displaystyle\inte \psi_{1}^2\not=0,
\]
and moreover,
\[
-\displaystyle\frac{2+p}{2} \beta _k+\alp _k \Big(2\displaystyle\inte w\psi_{3}+\displaystyle\inte \psi_{1}^2\Big)=0,\ \ i.e.,\ \beta _k=C^*\alp_k,
\]
where $C^*\not =0$ is as in (\ref{3:35a}). Finally, the expansion (\ref{3:36}) follows directly from (\ref{3:35}) and Lemma \ref{lem3.2}, and we are done.
\qed


We remark from (\ref{2:beta}) and Proposition \ref{prop3.4} that the  Lagrange multiplier $\mu _k\in \R$ of the
Euler-Lagrange equation (\ref{2:cond}) satisfies
\begin{equation}\label{3Aa:27}
\mu _k=-\frac{\lam}{\eps ^2_k}+\lam ^2C^*\eps_k^p +o(\eps_k^p)\,\ \mbox{as} \,\ k\to\infty,
\end{equation}
where $\lam>0$ is defined by (\ref{def:li}) with $g(0)=1$, and $C^*\not =0$ is given by (\ref{3:35a}). Moreover, following above results we finally conclude the following refined spike profiles.

\begin{thm}\label{thm3.5}
Suppose that $V(x)=h(x)\in C^2(\R^2)$ satisfies $\lim_{|x|\to\infty} h(x) = \infty$ and (\ref{1:H}) for  some $y_0\in\R^2$, where $h(x)$ is homogeneous of degree $p\ge 2$. If $u_a$ is a
positive minimizer of $e(a)$ for $a<a^*$. Then for any sequence $\{a_k\}$ with $a_k\nearrow a^*$ as $k\to\infty$, there exist a subsequence, still denoted by $\{a_k\}$, of $\{a_k\}$ and $\{x_k\}\subset\R^2$ such that the subsequence solution $ u_k=u_{a_k}$ satisfies for $\eps _k:=(a^*-a_k)^{\frac{1}{2+p}}$,
\begin{equation}\arraycolsep=1.5pt\begin{array}{lll}
u_k(x)&=& \displaystyle\frac { \lambda }{\|w\|_2} \Big\{  \displaystyle\frac{1}{\eps _k} w\Big(\frac{\lam (x-x_k)}{\eps _k}\Big)+\displaystyle\eps _k^{1+p}\Big[\psi_{1}+ C^*\psi_{2}\Big]\Big(\frac{\lam (x-x_k)}{\eps _k}\Big)\\[4mm]
&&+\displaystyle\eps _k^{3+2p}\Big[\psi_{3}+(C^*)^2\psi_{4} + C^* \psi_{5}\Big]\Big(\frac{\lam (x-x_k)}{\eps _k}\Big) \Big\} +o(\eps _k^{3+2p})    \,\ \mbox{as} \,\ k\to\infty \label{thm3.4:cona}
\end{array}\end{equation}
uniformly in $\R^2$,
where the unique maximum point $x_k$ of $u_k$ satisfies
\begin{equation}
 \Big|\frac{\lam x_k }{\eps _k }-y_0\Big|=\eps _k^{2+p}O( |y^0|)\,\ \mbox{as} \,\ k\to\infty
\label{thm3.4:conb}
\end{equation}
for some $y^0\in\R^2$, and $C^*\not =0$ is given by (\ref{3:35a}). Here $\psi_1(x),\cdots ,\psi_5(x)\in C^2(\R^2)\cap L^\infty(\R^2)$  are given in Lemma \ref{lem2.3} with $g(0)=1$ and Lemma \ref{lem3.2}.
\end{thm}

\noindent\textbf{Proof.} The refined spike profile (\ref{thm3.4:cona}) follows immediately from (\ref{2:cong}) and (\ref{3:36}). Also, Lemma \ref{lem3.1} and (\ref{3:35}) yield that the estimate (\ref{thm3.4:conb}) holds.
\qed

\vskip 0.07truein
\noindent\textbf{Proof of  Theorem \ref{1:thmB}.} Since the local uniqueness of Theorem \ref{1:thmA} implies that Theorem \ref{thm3.5} holds for the whole sequence $\{a_k\}$,  Theorem \ref{1:thmB} is proved.\qed

%


\section{Refined Spike Profiles: $V(x)=g(x)h(x)$}
The main purpose of  this section is to derive Theorem \ref{thm4.5} which extends the refined  spike behavior of Theorem \ref{1:thmB} to more general potentials $V(x)=g(x)h(x)\in C^2(\R^2)$, where $V(x)$ satisfies $\lim_{|x|\to\infty} V(x) = \infty$ and
\begin{enumerate}
\item[$(V).$]   $h(-x)=h(x)$ satisfies  (\ref{1:H}) and  is homogeneous of degree $p\ge 2$, $g(x)\in  C^m(\R^2)$ for some $2\le m\in \mathbb{N}\cup \{+\infty\}$ satisfies $0<C\le g(x)\le \frac{1}{C}$ in $\R^2$ and $ G(x):=g(x)-g(0)$,
    \[
 \mathcal{D} ^{\alp} G(0)=0\, \ \mbox{for all }\ |\alp|\le m-1, \ \mbox{and} \,\ \mathcal{D} ^\alp G(0)\not=0\,\ \mbox{for some }\ |\alp|=m.
    \]
\end{enumerate}
Here it takes $m=+\infty$ if $g(x)\equiv 1$.

\begin{rem} The property $h(-x)=h(x)$ in the above assumption $(V)$ implies that $y_0=0$ must occur in (\ref{1:H}).
\label{4:rem}\end{rem}

For the above type of potentials $V(x)$, suppose $\{u_k\}$ is a
positive minimizer sequence of $e(a_k)$ with $a_k\nearrow a^*$ as $k\to\infty$,  and let $w_k$ be defined by (\ref{2:cong}), where $x_k$ is the unique maximum point of $u_k$. Then Lemma \ref{lem2.3} still holds in this case, where $\alp _k>0$ and $\beta_k>0$ are defined in (\ref{2:beta}). Similar to Lemma \ref{lem3.1}, we start with the following estimates.

\begin{lem}\label{lem4.1} Suppose $V(x)=g(x)h(x)\in C^2(\R^2)$ satisfies $\lim_{|x|\to\infty} V(x) = \infty$ and the assumption $(V)$ for $p\ge 2$ and $2\le m\in \mathbb{N}\cup \{+\infty\}$. Then the unique maximum point $x_k$ of $u_k$ satisfies the following estimates:

\begin{enumerate}
\item If  $m$ is even, then we have
\begin{equation}\label{4:2bb}
\frac{\lam \alp _k|x_k|}{\eps _k}=o\big([\alp _k+\beta _k]^2+\alp_k\eps _k^m\big)\,\ \mbox{as} \,\ k\to\infty.
\end{equation}


\item If  $m$ is odd, then we have
\begin{equation}\label{4:2a}
\frac{\lam \alp _k|x_k|}{\eps _k}=O(\alp _k\eps _k^m|x_0|) +o\big([\alp _k+\beta _k]^2+\alp_k\eps _k^m \big)\,\ \mbox{as} \,\ k\to\infty,
\end{equation}
where $x_0\in \R^2$ satisfies
\begin{equation}\label{4:2b}
g(0)\inte\frac{\partial w}{\partial x_1} \big[x_0\cdot \nabla h(x)\big]w+\frac{1}{\lam ^m}\sum _{|\alp|=m} \inte\frac{\partial w}{\partial x_1}\Big[\frac{x^\alp}{\alp !}\mathcal{D}^{\alp}g(0)\Big]h(x)w=0.
\end{equation}

\end{enumerate}
\end{lem}

\noindent\textbf{Proof.} Recall that $\psi_{1}(x)$ and $\psi_{2}(x)$ are given in Lemma \ref{lem2.3}.  Since $h(-x)=h(x)$, we have $\psi_i(-x)=\psi_i(x)$ for $i=1,\,2$ and thus
\begin{equation}\label{H:3}
\inte\frac{\partial w}{\partial x_1}\psi_1 =\displaystyle\inte\frac{\partial w}{\partial x_1}w\psi _1^2=\displaystyle\inte\frac{\partial w}{\partial x_1}w\psi _1\psi _2=0.
\end{equation}
Since (\ref{1:H}) holds with $y_0=0$ as shown in Remark \ref{4:rem}, the same calculations of (\ref{3a:21})--(\ref{3a:23}) then yield that
\begin{equation}\label{4:3}\arraycolsep=1.5pt\begin{array}{lll}
&& o(\alp _k^2+\alp_k\beta _k)=\displaystyle\inte\frac{\partial w}{\partial x_1}\mathcal{L}_kw_k\\[4mm]
&=& o(\alp _k\beta _k+\beta^2 _k)-\alp _k\displaystyle\inte\frac{\partial w}{\partial x_1}\Big[\frac{1}{a^*}\bar u_k^3 +\displaystyle\frac{\bar u_k }{\lam ^{2+p}} g\big(\frac{\eps_kx}{\lam}+x_k\big)h\big(x+\frac{\lam x_k}{\eps _k}\big)\Big]\\[4mm]
&=& o(\alp _k\beta _k+\beta^2 _k)-\displaystyle\frac{\alp _k}{a^*}\inte\frac{\partial w}{\partial x_1}\big( \bar u_k^3-w^3\big) \\[4mm]
&& -\displaystyle\frac{\alp _k}{\lam ^{2+p}}\inte \frac{\partial w}{\partial x_1}\Big[g\big(\frac{\eps_kx}{\lam}+x_k\big)h\big(x+\frac{\lam x_k}{\eps _k}\big)\bar u_k -g(0)h(x)w\Big]\\[4mm]
&=& o(\alp _k\beta _k+\beta^2 _k)-\displaystyle\frac{\alp _k}{\lam ^{2+p}}\inte \frac{\partial w}{\partial x_1}\Big[g\big(\frac{\eps_kx}{\lam}+x_k\big)h\big(x+\frac{\lam x_k}{\eps _k}\big)\bar u_k -g(0)h(x)w\Big]\\[4mm]
&=& o(\alp _k\beta _k+\beta^2 _k)-I_1,
\end{array}\end{equation}
where the first equality follows from (\ref{3a:22}) and (\ref{H:3}). Similar to (\ref{3a:24}), we deduce from (\ref{1:H}) with $y_0=0$ that
\[\arraycolsep=1.5pt\begin{array}{lll}
\displaystyle\frac{\lam ^{2+p}}{\alp _k} I_1&=&\displaystyle\inte\frac{\partial w}{\partial x_1}\Big\{g(0)h\big(x+\frac{\lam x_k}{\eps _k}\big)\big[\bar u_k-w\big]+g(0)\big[h\big(x+\frac{\lam x_k}{\eps _k}\big)-h(x)\big]w\Big\}\\[4mm]
&&+\displaystyle\inte\frac{\partial w}{\partial x_1}\Big[g\big(\frac{\eps_kx}{\lam}+x_k\big)- g(0)\Big]h\big(x+\frac{\lam x_k}{\eps _k}\big)\bar u_k\\[4mm]
&=&o(\alp _k+\big|\frac{\lam x_k}{\eps _k}\big|+\beta _k)+g(0)\displaystyle\inte\frac{\partial w}{\partial x_1}\Big(\frac{\lam x_k}{\eps _k}\cdot \nabla h(x)\Big)w\\[4mm]
&&+\displaystyle\Big(\frac{\eps_k}{\lam}\Big)^m \sum _{|\alp|=m}  \inte\frac{\partial w}{\partial x_1}\Big[\frac{1}{\alp !}\Big(x+\frac{\lam x_k}{\eps _k}\Big)^\alp \mathcal{D}^{\alp}g(0)\Big]h\big(x+\frac{\lam x_k}{\eps _k}\big)\bar u_k+o(\eps_k^m),
\end{array}\]
which then implies that
\begin{equation}\label{4:7}\arraycolsep=1.5pt\begin{array}{lll}
I_1&=&\displaystyle\frac{\alp _k}{\lam ^{2+p}}\Big(\frac{\eps_k}{\lam}\Big)^m\sum _{|\alp|=m}\inte\frac{\partial w}{\partial x_1}\Big[\frac{1}{\alp !}\Big(x+\frac{\lam x_k}{\eps _k}\Big)^\alp \mathcal{D} ^{\alp}g(0)\Big]h\big(x+\frac{\lam x_k}{\eps _k}\big)\bar u_k\\[4mm]
&&+\displaystyle\frac{\alp _k}{\lam ^{2+p}}g(0)\inte\frac{\partial w}{\partial x_1} \Big(\frac{\lam x_k}{\eps _k}\cdot \nabla h(x)\Big)w+o(\alp _k^2+\alp_k\beta _k+\big|\frac{\lam x_k}{\eps _k}\big|+\alp_k\eps_k^m).
\end{array}\end{equation}
Combining (\ref{4:3}) and (\ref{4:7}), we then conclude from the estimate (\ref{2N:a3})  that
\begin{equation}\label{4:8}
\arraycolsep=1.5pt\begin{array}{lll}
&& \displaystyle\frac{\alp _k}{\lam ^{2+p}}g(0)\inte\frac{\partial w}{\partial x_1} \Big(\frac{\lam x_k}{\eps _k}\cdot \nabla h(x)\Big)w\\[4mm]
&=&-\displaystyle\frac{\alp _k}{\lam ^{2+p}}\Big(\frac{\eps_k}{\lam}\Big)^m\sum _{|\alp|=m}\inte\frac{\partial w}{\partial x_1}\Big[\frac{x^\alp}{\alp !}\mathcal{D}^{\alp}g(0)\Big]h(x)w+o\Big([\alp _k+\beta _k]^2+\alp _k\eps_k^m\Big).
\end{array}\end{equation}

If $m$ is even, one can note that
\[\sum _{|\alp|=m}\inte\frac{\partial w}{\partial x_1}\Big[\frac{x^\alp}{\alp !}\mathcal{D} ^{\alp}g(0)\Big] h(x)w=0,\]
and it then follows from (\ref{4:8}) and (\ref{1:H}) with $y_0=0$ that (\ref{4:2bb}) holds. If $m$ is odd, we then derive from  (\ref{4:8}) that both (\ref{4:2a}) and (\ref{4:2b}) hold.
\qed

\begin{lem}\label{lem4.2} Suppose $V(x)=g(x)h(x)\in C^2(\R^2)$ satisfies $\lim_{|x|\to\infty} V(x) = \infty$ and the assumption $(V)$ for $p\ge 2$ and $2\le m\in \mathbb{N}\cup \{+\infty\}$. Let $\psi_1(x)$ and $\psi_2(x)$ be given in Lemma \ref{lem2.3} with $y_0=0$.
Then $w_k$ satisfies
 \begin{equation}\label{lem4.2:b1}\arraycolsep=1.5pt\begin{array}{lll}
w_k :&=&\alp _k\psi_{1} +\beta _k\psi_{2} + \alp ^2_k\psi_{3} +\beta ^2_k\psi_{4}\\[3mm]
&&+\alp _k \eps_k^m\phi +\alp _k\beta _k\psi_{5}+o\big([\alp_k+\beta _k]^2+\alp _k\eps _k^m\big) \ \ \mbox{as} \ \ k\to\infty,
\end{array}\end{equation}
where $\psi_{i}(x)\in C^2(\R^2)\cap L^\infty(\R^2)$, $i=3,4,5$, solves uniquely
     \begin{equation}\label{lem4.2:a2}
 \nabla \psi_{i}(0)=0\ \, \mbox{and}\,\ \mathcal{L}\psi_{i}(x)=g_i(x)
 \,\ \mbox{in}\,\ \R^2,\,\ i=3,4,5,
   \end{equation}
 and $g_i(x)$ satisfies
 \begin{equation}\label{lem4.2:a3}\arraycolsep=1.5pt
 g_i(x)=\left\{\begin{array}{lll}
 &3w\psi^2_1-\displaystyle\Big(\frac{3w^2}{a^*}+\frac{g(0)h(x)}{\lam ^{2+p}}\Big)\psi_1,\,\ &\mbox{if}\ \ i=3;\\[3mm]
 &3w\psi^2_2+\psi_2, \,\ &\mbox{if}\ \ i=4;\\[2mm]
 &6w\psi_1\psi_2+\psi _1-\displaystyle\Big(\frac{3w^2}{a^*}+\frac{g(0)h(x)}{\lam ^{2+p}}\Big)\psi_2, \,\ &\mbox{if}\ \ i=5.
 \end{array}\right.
    \end{equation}
Here $\phi\in C^2(\R^2)\cap L^\infty(\R^2)$ solves uniquely
     \begin{equation}\label{lem4.2:b2}\arraycolsep=1.5pt\begin{array}{lll}
\mathcal{L}\phi (x)&=&-\displaystyle\frac{1}{\lam ^{2+p}}\Big\{\big[x_0\cdot \nabla h(x)\big]g(0)w \\[3mm]
&&\qquad\quad \ \, + \displaystyle\frac{1}{\lam ^m}   \sum _{|\alp|=m}  \Big[\frac{x^\alp}{\alp !}\mathcal{D} ^{\alp}g(0)\Big]h(x)w\Big\}
 \,\ \mbox{in}\,\ \R^2, \ \, \mbox{and}\,\  \nabla \phi (0)=0,
\end{array}   \end{equation}
where $x_0=0$ holds for the case where $m$ is even, and $x_0\in \R^2$ satisfies (\ref{4:2b}) for the case where $m$ is odd.
\end{lem}

\noindent\textbf{Proof.} Following Lemma \ref{lem2.3}(3), we set
\[
v_k=w_k-\alp_k\psi_1-\beta_k\psi_2.
\]
Similar to (\ref{3a:31a}), we then have
     \begin{equation}\label{lem4.2:b3}\arraycolsep=1.5pt\begin{array}{lll}
 \mathcal{L}_kv_k&=&w_k(\alp_k\psi_1+\beta_k\psi_2)(3w+w_k)+\beta_kw_k-  \displaystyle\frac{\alp_k}{a^*}(\bar u_k^3-w^3)\\[1mm]
&&- \displaystyle\frac{\alp_k}{\lam ^{2+p}}\Big[g\big(\frac{\eps_kx}{\lam}+x_k\big)h\big(x+\frac{\lam x_k}{\eps _k}\big)\bar u_k -g(0)h(x)w\Big]\\[3mm]
&=&w_k(\alp_k\psi_1+\beta_k\psi_2)(3w+w_k)+\beta_kw_k\\[2mm]
&&-\displaystyle\frac{\alp_k}{a^*}w_k(3w^2+3ww_k+w_k^2)-I_2,
\end{array} \end{equation}
where $I_2$ satisfies
\begin{equation}\label{lem4.2:b4}\arraycolsep=1.5pt\begin{array}{lll}
 I_2&=&\displaystyle\frac{\alp_k}{\lam ^{2+p}}\Big\{\Big[g\big(\frac{\eps_kx}{\lam}+x_k\big)-g(0)\Big]h\big(x+\frac{\lam x_k}{\eps _k}\big)\bar u_k \\[3mm]
 &&\qquad\quad +\displaystyle g(0)\Big[h\big(x+\frac{\lam x_k}{\eps _k}\big)-h(x)\Big]\bar u_k +g(0)h(x)\big(\bar u_k-w\big)\Big\}\\[3mm]
&=& \displaystyle\frac{\alp_k}{\lam ^{2+p}}\Big\{\Big(\frac{\eps _k}{\lam}\Big)^m\sum _{|\alp|=m} \Big[\frac{1}{\alp !}\Big(x+\frac{\lam x_k}{\eps _k}\Big)^\alp \mathcal{D} ^{\alp}g(0)\Big]h\big(x+\frac{\lam x_k}{\eps _k}\big)\bar u_k\\[3mm]
&&\qquad\quad +\displaystyle g(0)\Big(\frac{\lam x_k}{\eps _k}\cdot \nabla h(x)\Big)\bar u_k +g(0)h(x)w_k\Big\}+\displaystyle o\Big(\frac{\alp _k x_k}{\eps _k}+\alp _k\eps _k^m\Big)\\[3mm]
&=&\alp_kw_k\displaystyle\frac{g(0)h(x)}{\lam ^{2+p}}
+\displaystyle\frac{\alp_k}{\lam ^{2+p}}\Big(\frac{\lam x_k}{\eps _k}\cdot \nabla h(x)\Big)g(0)\bar u_k \\[4mm]
&&+\displaystyle\frac{\alp_k\eps _k^m}{\lam ^{2+p+m}}\sum _{|\alp|=m} \Big[\frac{1}{\alp !}\Big(x+\frac{\lam x_k}{\eps _k}\Big)^\alp \mathcal{D}^{\alp}g(0)\Big]h\big(x+\frac{\lam x_k}{\eps _k}\big)\bar u_k
+o\big(\alp _k\eps _k^m\big),
\end{array} \end{equation}
where Lemma \ref{lem4.1} is used in the last equality. Applying  Lemma \ref{lem4.1} again, we then obtain from (\ref{lem4.2:b3}) and (\ref{lem4.2:b4}) that
\begin{equation}\label{4a:33a}\arraycolsep=1.5pt\begin{array}{lll}
 \mathcal{L}_kv_k&=&w_k(\alp_k\psi_1+\beta_k\psi_2)(3w+w_k)+\beta_kw_k-\displaystyle\frac{\alp_k}{a^*}w_k^2(3w+w_k)\\[2mm]
 &&- \alp_kw_k\displaystyle\Big[\frac{3w^2}{a^*}+\frac{g(0)h(x)}{\lam ^{2+p}}\Big]  -\displaystyle\frac{\alp_k}{\lam ^{2+p}}\Big(\frac{\lam x_k}{\eps _k}\cdot \nabla h(x)\Big)g(0)\bar u_k  \\[3mm]
 &&-\displaystyle\frac{\alp_k\eps _k^m}{\lam ^{2+p+m}}\sum _{|\alp|=m} \Big[\frac{1}{\alp !}\Big(x+\frac{\lam x_k}{\eps _k}\Big)^\alp \mathcal{D} ^{\alp}g(0)\Big]h\big(x+\frac{\lam x_k}{\eps _k}\big)\bar u_k+o\big(\alp _k\eps _k^m\big)\\[3mm]
 &=&\alp_k^2\displaystyle\Big[3w\psi^2_1-\Big(\frac{3w^2}{a^*}+\frac{g(0)h(x)}{\lam ^{2+p}}\Big)\psi_1 \Big]\\[3mm]
 &&+\alp_k\beta_k\displaystyle\Big[ 6w\psi_1\psi_2+\psi _1-\Big(\frac{3w^2}{a^*}+\frac{g(0)h(x)}{\lam ^{2+p}}\Big)\psi_2 \Big]\\[3mm]
 &&  -\displaystyle\frac{\alp_k\eps _k^m}{\lam ^{2+p }}\Big\{\big[x_0\cdot \nabla h(x)\big]g(0)w  + \displaystyle\frac{1}{\lam ^m}  \sum _{|\alp|=m}     \Big[\frac{x^\alp}{\alp !}\mathcal{D} ^{\alp}g(0)\Big]   h(x)w\Big\} \\[5mm]
 &&  +\beta_k^2(3w\psi^2_2+\psi_2) +o\big([\alp_k+\beta _k]^2+\alp _k\eps _k^m\big)\ \, \ \mbox{in}\ \ \R^2,
\end{array}
\end{equation}
where $x_0=0$ holds for the case where $m$ is even, and $x_0\in \R^2$ satisfies (\ref{4:2b}) for the case where $m$ is odd. Following (\ref{4a:33a}), the same argument of proving Lemma \ref{lem2.3} then gives (\ref{lem4.2:b1}), and the proof is therefore complete.
\qed

\begin{prop}\label{prop4.3} Suppose $V(x)=g(x)h(x)\in C^2(\R^2)$ satisfies $\lim_{|x|\to\infty} V(x) = \infty$ and the assumption $(V)$ for $p\ge 2$ and $2\le m\in \mathbb{N}\cup \{+\infty\}$. Let $\psi_1(x), \cdots, \psi_5(x)\in C^2(\R^2)\cap L^\infty(\R^2)$ be given in Lemma  \ref{lem2.3} with $y_0=0$ and Lemma \ref{lem4.2}, and $\phi$ is given by (\ref{lem4.2:b2}).
\begin{enumerate}
\item If  $ m>2+p$, then
\begin{equation}\label{4:35}
\beta _k=C^*\alp_k,
\end{equation}
and $w _k$ satisfies
\begin{equation}\label{4:36}
w_k :=\big[\psi_{1}+ C^*\psi_{2} \big] \alp_k+\big[\psi_{3}+(C^*)^2\psi_{4}+ C^* \psi_{5}\big]\alp_k^2+o(\alp_k^2) \,\ \mbox{as} \,\ k\to\infty,
\end{equation}
where the constant $C^*$ satisfies
\begin{equation}\label{4aa:37b}
C^*:=\frac{2}{2+p}\Big(2\displaystyle\inte w\psi_{3}+\displaystyle\inte \psi_{1}^2\Big)\not =0.
\end{equation}

\item If  $1\le m\le 2+p$  and $m$ is odd, then $\beta _k=C^*\alp_k$  and $w _k$ satisfies
\begin{equation}\label{4:36b}\arraycolsep=1.5pt\begin{array}{lll}
w_k :&=&\big[\psi_{1}+ C^*\psi_{2} \big] \alp_k+\phi\,\alp _k \eps_k^m\\[3mm]
&&+\big[\psi_{3}+(C^*)^2\psi_{4}+ C^* \psi_{5}\big]\alp_k^2+o(\alp _k \eps_k^m) \,\ \mbox{as} \,\ k\to\infty,
\end{array}\end{equation}
where the constant $C^*\not =0$ is given by (\ref{4aa:37b}).

\item  If  $1\le m<2+p$  and $m$ is even, consider
 \begin{equation}\label{4E:g}
\mathcal{S}=\sum _{|\alp|=m}   \inte \Big[\frac{x^\alp}{\alp !}\mathcal{D}  ^{\alp}g(0)\Big]h(x)w^2.
\end{equation}
Then for the case where $\mathcal{S}=0$, we have $\beta _k=C^*\alp_k$  and $w _k$ satisfies (\ref{4:36b}),  where the constant $C^*\not =0$ is given by (\ref{4aa:37b}). However, for the case where $\mathcal{S}\neq 0$, we have
 \begin{equation}\label{4:35c}
\beta _k=C_1^*\eps_k^m,
\end{equation}
and $w _k$ satisfies
\begin{equation}\label{4:36c}
w_k :=C_1^*\psi_{2} \eps_k^m  + \psi_{1}\alp_k+ (C_1^*)^2\psi_{4} \eps_k^{2m}+
o(\eps_k^{\min \{2+p,2m\}}) \,\ \mbox{as} \,\ k\to\infty,
 \end{equation}
where the constant $C_1^*$ satisfies
\begin{equation}\label{4:35cc}
C_1^*=-\displaystyle\frac{m+p}{(2+p)\lam ^{2+p+m}}\sum _{|\alp|=m}   \inte \Big[\frac{x^\alp}{\alp !}\mathcal{D}  ^{\alp}g(0)\Big]h(x)w^2\not =0.
\end{equation}

\item  If  $m=2+p$ is even, then
\begin{equation}\label{4:35d2}
\beta _k=C_2^*\alp_k,
\end{equation}
and $w _k$ satisfies
\begin{equation}\label{4:35d3}\arraycolsep=1.5pt\begin{array}{lll}
w_k :&=& \big[\psi_{1}+ C_2^*\psi_{2} \big] \alp_k+\big[\psi_{3}+(C_2^*)^2\psi_{4}+ C_2^* \psi_{5}+\phi\big]\alp_k^2+o(\alp_k^2) \,\ \mbox{as} \,\ k\to\infty,
\end{array}\end{equation}
where the constant $C_2^*$ satisfies
\begin{equation}\label{4:35d1}
C_2^*= \frac{2}{2+p}\Big[2\displaystyle\inte w\psi_{3}+\displaystyle\inte \psi_{1}^2+2\inte  w\phi\Big]\not =0.
\end{equation}
\end{enumerate}
\end{prop}

\noindent\textbf{Proof.} The same argument of proving Lemma  \ref{lem3.3} with $y_0=0$ yields that
\begin{equation}\label{4a:37a}
 \displaystyle\inte w\psi_{1} =0,\ \ \inte w\psi_{2} =0 \,\ \mbox{and}\,\  I= \inte \big(2w\psi_4+\psi^2_2\big)=0,
\end{equation} and
 \begin{equation}\label{4a:claim3}
II=2\inte w\psi_5+2\inte \psi_1\psi_2=-\frac{2+p}{2}<0.
\end{equation}
It thus follows from (\ref{2a:k}) and Lemma \ref{lem4.2} that
\begin{equation}\label{4aa:37}\arraycolsep=1.5pt\begin{array}{lll}
0&=&2\displaystyle\inte ww_k+\displaystyle\inte w_k^2\\[3mm]
&=&2\displaystyle\inte w\big(\alp _k\psi_{1}  +\beta _k\psi_{2}+ \alp ^2_k\psi_{3} +\beta ^2_k\psi_{4}+\alp _k \eps_k^m\phi +\alp _k\beta _k\psi_{5}\big)\\[4mm]
&&+\displaystyle\inte \big(\alp _k\psi_{1} +\beta _k\psi_{2} + \alp ^2_k\psi_{3} +\beta ^2_k\psi_{4}+\alp _k \eps_k^m\phi +\alp _k\beta _k\psi_{5}\big)^2\\[4mm]
&&+o\big([\alp_k+\beta _k]^2+\alp _k\eps _k^m\big)\\[4mm]
&=& \alp _k\Big(\displaystyle2\inte  w\psi_{1}\Big)+\beta _k\Big(\displaystyle 2\inte  w\psi_{2}\Big)+\beta _k^2\Big(2\displaystyle\inte  w\psi_{4}+\displaystyle\inte \psi_{2}^2\Big)\\[4mm]
&&+\alp_k\beta _k\Big(2\displaystyle\inte w\psi_{5}+2\displaystyle\inte \psi_{1}\psi_{2}\Big)+\alp _k^2\Big(2\displaystyle\inte w\psi_{3}+\displaystyle\inte \psi_{1}^2\Big)\\[4mm]
&&+\alp _k\eps_k^m\Big(\displaystyle2\inte  w\phi\Big)+o\big([\alp_k+\beta _k]^2+\alp _k\eps _k^m\big)\\[4mm]
&=&-\displaystyle\frac{2+p}{2}\alp_k\beta _k+\alp _k^2\Big(2\displaystyle\inte w\psi_{3}+\displaystyle\inte \psi_{1}^2\Big)+\alp _k\eps_k^m\Big(\displaystyle2\inte  w\phi\Big)\\[4mm]
&&+o\big([\alp_k+\beta _k]^2+\alp _k\eps _k^m\big),
\end{array}
\end{equation}
where (\ref{4a:37a}) and (\ref{4a:claim3}) are used in the last equality. Following (\ref{4aa:37}), we next carry out the proof by considering separately the following four cases:
\vskip 0.1truein

\noindent{\em  Case 1. $m>2+p$.} In this case, it follows from (\ref{4aa:37}) that the constant $C^*$ defined in (\ref{4aa:37b}) is nonzero and
\[
-\displaystyle\frac{2+p}{2} \beta _k+\alp _k \Big(2\displaystyle\inte w\psi_{3}+\displaystyle\inte \psi_{1}^2\Big)=0,\ \ i.e.,\ \beta _k=C^*\alp_k.
\]
Moreover, the expansion (\ref{4:36}) follows directly from (\ref{4:35}) and Lemma \ref{lem4.2}, and Case 1 is therefore proved.

\vskip 0.1truein

\noindent{\em  Case 2. $1\le m\le 2+p$ and $m$ is odd.}  In this case, since $m$ is odd and $h(-x)=h(x)$, we obtain from (\ref{lem2.1:4a}) and (\ref{lem4.2:b2}) that
\[\arraycolsep=1.5pt\begin{array}{lll}
2\displaystyle\inte  w\phi &=&2\displaystyle\inte \phi\mathcal{L}\psi_{2}=2\displaystyle\inte \psi_{2}\mathcal{L}\phi\\[4mm]
&=&\displaystyle\frac{1}{\lam ^{2+p}}\inte \Big\{\big[x_0\cdot \nabla h(x)\big]g(0)w \\[3mm]
&&+ \displaystyle\frac{1}{\lam ^m}  \sum _{|\alp|=m}  \Big[\frac{x^\alp}{\alp !}\mathcal{D} ^{\alp}g(0)\Big]  h(x)w\Big\}\big(w+x\cdot\nabla w\big)=0.
\end{array}\]
We then derive from (\ref{4aa:37}) that (\ref{4aa:37b}) still holds and thus $\beta _k=C^*\alp_k$. Further, the expansion (\ref{4:36b}) follows directly from (\ref{lem4.2:b1}) and (\ref{4:35}).

\vskip 0.1truein
\noindent{\em  Case 3. $1\le m<2+p$ and $m$ is even.}  Since $m$ is even, then $x_0=0$ holds in (\ref{lem4.2:b2}). Further, since $x^\alp h(x)$ is homogeneous of degree $m+p$, we then obtain from (\ref{lem4.2:b2}) that
\begin{equation}\label{40:21}\arraycolsep=1.5pt\begin{array}{lll}
2\displaystyle\inte  w\phi &=&2\displaystyle\inte \phi\mathcal{L}\psi_{2}=2\displaystyle\inte \psi_{2}\mathcal{L}\phi\\[4mm]
&=&\displaystyle\frac{1}{\lam ^{2+p+m}}\sum _{|\alp|=m}   \inte \Big[\frac{x^\alp}{\alp !}\mathcal{D }^{\alp}g(0)\Big]h(x)w \big(w+x\cdot\nabla w\big) \\[4mm]
&=&\displaystyle\frac{1}{\lam ^{2+p+m}}\sum _{|\alp|=m}   \inte \Big[\frac{x^\alp}{\alp !}\mathcal{D} ^{\alp}g(0)\Big]h(x)w^2\\[4mm]
&&
+\displaystyle\frac{1}{2\lam ^{2+p+m}}\sum _{|\alp|=m}   \inte \Big[\frac{x^\alp}{\alp !}\mathcal{D} ^{\alp}g(0)\Big]h(x)\big(x\cdot\nabla w^2 \big)\\[4mm]
&=&\displaystyle\frac{1}{\lam ^{2+p+m}}\sum _{|\alp|=m}   \inte \Big[\frac{x^\alp}{\alp !}\mathcal{D} ^{\alp}g(0)\Big]h(x)w^2\\[4mm]
&&
-\displaystyle\frac{1}{2\lam ^{2+p+m}}\sum _{|\alp|=m}   \inte w^2\Big\{2\Big[\frac{x^\alp}{\alp !}\mathcal{D}  ^{\alp}g(0) h(x)\Big]\\[4mm]
&&+\displaystyle x\cdot\nabla \Big[\frac{x^\alp}{\alp !}\mathcal{D}  ^{\alp}g(0) h(x)\Big]\Big\}\\[4mm]
&=&
-\displaystyle\frac{m+p}{2\lam ^{2+p+m}}\sum _{|\alp|=m}   \inte \Big[\frac{x^\alp}{\alp !}\mathcal{D}  ^{\alp}g(0)\Big]h(x)w^2:=-\displaystyle\frac{m+p}{2\lam ^{2+p+m}}\mathcal{S},
\end{array}\end{equation}
where $\mathcal{S}$ is as in (\ref{4E:g}). Therefore, if $\mathcal{S}=0$, then we are in the same situation as that of above Case 2, which gives that  $\beta _k=C^*\alp_k$ and $w _k$ satisfies (\ref{4:36b}), where the constant $C^*\not =0$ is given by (\ref{4aa:37b}).

We next consider the case where $\mathcal{S}\neq 0$. By applying (\ref{40:21}), in this case we derive from (\ref{4aa:37}) that
\[-\displaystyle\frac{2+p}{2}\alp_k\beta _k+\alp _k\eps_k^m\Big(\displaystyle2\inte  w\phi\Big)=0,
\]
which implies that $\beta _k=C_1^*\eps_k^m$, where the constant $C_1^*\not =0$ satisfies (\ref{4:35cc}) in view of (\ref{40:21}).
Further, the expansion (\ref{4:36c}) follows directly from (\ref{4:35c}) and Lemma \ref{lem4.2}.

\vskip 0.1truein
\noindent{\em  Case 4. $m=2+p$  is even.}  In this case, we derive from (\ref{4aa:37}) that
\[-\displaystyle\frac{2+p}{2}\alp_k\beta _k+\alp _k^2\Big(2\displaystyle\inte w\psi_{3}+\displaystyle\inte \psi_{1}^2+\displaystyle2\inte  w\phi\Big)=0,
\]
which gives that $\beta _k=C_2^*\alp_k $, where the constant $C_2^*\not =0$ satisfies (\ref{4:35d1}).
Further, the expansion (\ref{4:35d3}) follows directly from (\ref{4:35d2}) and Lemma \ref{lem4.2}.
\qed


Applying directly Lemmas \ref{lem4.1} and \ref{lem4.2} as well as Proposition \ref{prop4.3}, we now conclude the following main results of this section. Recall that   $\lam>0$ is defined by (\ref{def:li}) with $y_0=0$, $\psi_1(x), \cdots, \psi_5(x)\in C^2(\R^2)\cap L^\infty(\R^2)$ are given in Lemma  \ref{lem2.3} with $y_0=0$ and Lemma \ref{lem4.2}, and $\phi$ is given by (\ref{lem4.2:b2}).

\begin{thm}\label{thm4.5}
Suppose $V(x)=g(x)h(x)\in C^2(\R^2)$ satisfies $\lim_{|x|\to\infty} V(x) = \infty$ and the assumption $(V)$ for $p\ge 2$ and $2\le m\in \mathbb{N}\cup \{+\infty\}$.  Let $u_a$ be a positive minimizer of (\ref{f}) for $a<a^*$. Then for any sequence $\{a_k\}$ with $a_k\nearrow a^*$ as $k\to\infty$, there exists a subsequence, still denoted by $\{a_k\}$, of $\{a_k\}$ such that $u_k=u_{a_k}$ has a unique maximum point $x_k\in\R^2$ and satisfies for $\eps _k:=(a^*-a_k)^{\frac{1}{2+p}}$,
\begin{enumerate}
\item If $ m>2+p$, then we have
\begin{equation}\arraycolsep=1.5pt\begin{array}{lll}
u_k(x)&=& \displaystyle\frac { \lambda }{\|w\|_2} \Big\{  \displaystyle\frac{1}{\eps _k} w\Big(\frac{\lam (x-x_k)}{\eps _k}\Big)+\displaystyle\eps _k^{1+p}\Big[\psi_1+C^*\psi_2\Big]\Big(\frac{\lam (x-x_k)}{\eps _k}\Big)\\[4mm]
&&+\displaystyle\eps _k^{3+2p}\Big[\psi_{3}+(C^*)^2\psi_{4}+ C^* \psi_{5}\Big]\Big(\frac{\lam (x-x_k)}{\eps _k}\Big)\Big\}  +o(\eps _k^{3+2p})      \,\ \mbox{as} \,\ k\to\infty
\end{array}\label{thm4.5:A}\end{equation}
uniformly in $\R^2$, where $x_k$  satisfies
\begin{equation}
\frac{|x_k|}{\eps _k}=O(\eps _k^m |y^0|)+o(\eps _k^{2+p})\,\ \mbox{as} \,\ k\to\infty
\label{thm4.5:B}
\end{equation}
for some $y^0\in\R^2$, and the constant $C^*\not =0$ is given by (\ref{4aa:37b}). Further, if $m$ is even, then  $x_k$ satisfies
\begin{equation}
\frac{|x_k|}{\eps _k^{3+p}}= o(1)\,\ \mbox{as} \,\ k\to\infty.
\label{thm4.5:Y}
\end{equation}

\item If  $1\le m\le 2+p$  and $m$ is odd, then we have
\begin{equation}\arraycolsep=1.5pt\begin{array}{lll}
&u_k(x)= \displaystyle\frac { \lambda }{\|w\|_2} \Big\{& \displaystyle\frac{1}{\eps _k} w\Big(\frac{\lam (x-x_k)}{\eps _k}\Big)+\displaystyle\eps _k^{1+p}\Big[\psi_1+C^*\psi_2\Big]\Big(\frac{\lam (x-x_k)}{\eps _k}\Big)\\[4mm]
&&+\displaystyle\eps _k^{3+2p}\big[\psi_{3}+(C^*)^2\psi_{4}+ C^* \psi_{5}\big]\Big(\frac{\lam (x-x_k)}{\eps _k}\Big)\\[4mm]
&&+\displaystyle\eps _k^{1+m+p}\phi\Big(\frac{\lam (x-x_k)}{\eps _k}\Big)\Big\}   +o(\eps _k^{1+m+p})      \,\ \mbox{as} \,\ k\to\infty
\end{array}\label{thm4.5:C}\end{equation}
uniformly in $\R^2$, where $x_k$  satisfies
\begin{equation}
\frac{|x_k|}{\eps _k^{m+1}}= O(|y^0|)\,\ \mbox{as} \,\ k\to\infty.
\label{thm4.5:YY}
\end{equation}
for some $y^0\in\R^2$, and the constant $C^*\not =0$ is given by (\ref{4aa:37b}).

\item  If  $m=2+p$ is even, then we have
\begin{equation}\arraycolsep=1.5pt\begin{array}{lll}
u_k(x)&=& \displaystyle\frac { \lambda }{\|w\|_2} \Big\{  \displaystyle\frac{1}{\eps _k} w\Big(\frac{\lam (x-x_k)}{\eps _k}\Big)+\displaystyle\eps _k^{1+p}\Big[\psi_1+C^*_2\psi_2\Big]\Big(\frac{\lam (x-x_k)}{\eps _k}\Big)\\[4mm]
&&+\displaystyle\eps _k^{3+2p}\Big[\psi_{3}+(C_2^*)^2\psi_{4}+ C_2^* \psi_{5}+\phi\Big]\Big(\frac{\lam (x-x_k)}{\eps _k}\Big)\Big\}  +o(\eps _k^{3+2p})    \,\ \mbox{as} \,\ k\to\infty
\end{array}\label{thm4.5:E}\end{equation}
uniformly in $\R^2$, where $x_k$  satisfies (\ref{thm4.5:Y}) and the constant $C_2^*\not =0$ is defined by (\ref{4:35d1}).

\item If  $1\le m<2+p$  and $m$ is even, let the constant $\mathcal{S}$ be defined in (\ref{4E:g}).
Then for the case where $\mathcal{S}=0$, $u_k$ satisfies (\ref{thm4.5:C}) and $x_k$  satisfies
 \begin{equation}
\frac{|x_k|}{\eps _k^{m+1}}= o(1)\,\ \mbox{as} \,\ k\to\infty.
\label{thm4.5:YYY}
\end{equation}
However, for the case where $\mathcal{S}\neq 0$, $u_k$ satisfies
\begin{equation}\arraycolsep=1.5pt\begin{array}{lll}
u_k(x)&=& \displaystyle\frac { \lambda }{\|w\|_2} \Big\{  \displaystyle\frac{1}{\eps _k} w\Big(\frac{\lam (x-x_k)}{\eps _k}\Big)+\displaystyle\eps _k^{m-1}C^*_1\psi_2 \Big(\frac{\lam (x-x_k)}{\eps _k}\Big)\\[4mm]
&&+\displaystyle\eps _k^{2m-1}(C^*_1)^2\psi_4 \Big(\frac{\lam (x-x_k)}{\eps _k}\Big)+\displaystyle\eps _k^{1+p}\psi_1\Big(\frac{\lam (x-x_k)}{\eps _k}\Big)\Big\} \\[4mm]
&& +o(\eps_k^{\min \{2+p,2m\}-1})    \,\ \mbox{as} \,\ k\to\infty
\end{array}\label{thm4.5:F}\end{equation}
uniformly in $\R^2$, where  $x_k$  satisfies (\ref{thm4.5:YYY}), and the constant $C_1^*\not =0$ is defined by (\ref{4:35cc}).
\end{enumerate}
\end{thm}

\noindent\textbf{Proof.} (1). If $ m>2+p$, then (\ref{thm4.5:A}) follows directly from Proposition \ref{prop4.3}(1), and (\ref{thm4.5:B}) follows from Lemma \ref{lem4.1}. Specially, if $m$ is even, then Lemma \ref{lem4.1} gives $y^0=0$, and therefore (\ref{thm4.5:B}) implies (\ref{thm4.5:Y}).

(2).\, If $1\le m\le 2+p$  and $m$ is odd, then Proposition \ref{prop4.3}(2) gives (\ref{thm4.5:C}). Moreover, it yields from   (\ref{4:2a}) that  $x_k$ satisfies $\Big|\frac{x_k}{\eps _k}\Big|=O(\eps _k^m |y^0|)+o(\eps _k^m)$ as $k\to\infty$, which then implies (\ref{thm4.5:YY}) for some $y^0\in\R^2$.

(3).\, If  $m=2+p$ is even, then Proposition \ref{prop4.3}(4) gives (\ref{thm4.5:E}), and we reduce from (\ref{4:2bb}) that $x_k$ satisfies (\ref{thm4.5:Y}).

(4).\, If  $1\le m<2+p$  and $m$ is even, it then follows from (\ref{4:2bb}) that $x_k$ always satisfies (\ref{thm4.5:YYY}). Moreover, Proposition \ref{prop4.3}(3) gives that if $\mathcal{S}=0$, then $u_k$ satisfies (\ref{thm4.5:C}); if $\mathcal{S}\neq 0$, then $u_k$ satisfies (\ref{thm4.5:F}).
\qed


\appendix

\section{Appendix: The Proof of Lemma \ref{lem3.3}}
In this appendix, we shall follow Lemmas \ref{lem2.3} and  \ref{lem3.2} to address the proof of Lemma \ref{lem3.3}, i.e., (\ref{3a:37a})--(\ref{a:claim3}).\vskip 0.1truein

\noindent{\em The proof of (\ref{3a:37a}).} Under the assumptions of Lemma \ref{lem3.3}, we first note that the equation (\ref{lem2.1:2})  can be simplified as
\begin{equation}\label{00:21}
 \nabla \psi_{1}(0)=0,\quad \mathcal{L}\psi_{1}=- \frac{2w^3}{\inte w^4}-\frac{2h(x+y_0)w}{p\inte h(x+y_0)w^2}
 \,\ \mbox{in}\,\ \R^2,
\end{equation}
due to the fact that
\begin{equation}\label{00a:21}
a^*= \|w\|_2^2=\frac 1 2 \inte w^4.
\end{equation}
By (\ref{1:H}), (\ref{lem2.1:4a}) and (\ref{00:21}), we then have
 \[\arraycolsep=1.5pt\begin{array}{lll}
 2\displaystyle\inte w\psi _1&=&2\displaystyle\inte  \displaystyle\mathcal{L}\psi_{2}\psi _1=2\displaystyle\inte \psi_{2} \mathcal{L}\psi _1\\[4mm]
 &=&  \displaystyle\inte   \Big[\frac{2w^3}{\inte w^4}+\frac{2h(x+y_0)w}{p\inte h(x+y_0)w^2}\Big]\big(w+x\cdot \nabla w\big)\\[4mm]
 &=&2+\displaystyle\frac{2}{p}+\frac{2}{\inte w^4}\inte w^3\big(x\cdot \nabla w\big)+\frac{2 }{p\inte h(x+y_0)w^2}\inte h(x+y_0)w\big(x\cdot \nabla w\big)\\[4mm]
 &=&2+\displaystyle\frac{2}{p}+\frac{1  }{2\inte w^4}\inte \big(x\cdot \nabla w^4\big) +\frac{1}{p\inte h(x+y_0)w^2}\inte  h(x+y_0)\big(x\cdot \nabla w^2\big)\\[4mm]
 &=&2+\displaystyle\frac{2}{p}-1-\frac{1}{p\inte h(x+y_0)w^2} \inte  w^2\Big[2h(x+y_0)+\big(x\cdot \nabla h(x+y_0)\big)\Big]\\[4mm]
  &=&2+\displaystyle\frac{2}{p}-1-\frac{(p+2)}{p}=0,
 \end{array}\]
since $(x+y_0)\cdot \nabla h(x+y_0)=ph(x+y_0)$ and $\inte  w^2 \big[y_0\cdot \nabla h(x+y_0)\big]=0$. Also, we deduce from  (\ref{lem2.1:4a}) that
\[
2\inte w\psi_2=-\inte w(w+x\cdot \nabla w)=-\inte w^2-\frac 1 2 \inte  (x\cdot \nabla w^2)=0,
\]
which thus completes the proof of  (\ref{3a:37a}).
\qed

\vskip 0.1truein

\noindent{\em The proof of (\ref{3a:37aa}).} By Lemmas \ref{lem2.3} and  \ref{lem3.2}, we obtain that
\[\arraycolsep=1.5pt\begin{array}{lll}
  I&=&\displaystyle\inte \big(2w\psi_4+\psi^2_2\big)=\displaystyle\inte\psi^2_2+2\langle  \mathcal{L}\psi _2,\psi_4\rangle\\[4mm]
  &=&\displaystyle\inte\psi^2_2+2\langle  \psi _2,\mathcal{L}\psi_4\rangle
=\displaystyle\inte\psi^2_2+2\langle  \psi _2,(3w\psi_2^2+\psi_2)\rangle\\[4mm]
&=&3\displaystyle\inte\psi^2_2+\displaystyle6\inte w\psi_2^3,
  \end{array}\]
which implies that
\begin{equation}\label{0:21}
\arraycolsep=1.5pt\begin{array}{lll}
  \displaystyle\frac{4I}{3}\div 2\pi&=&4\Big[\displaystyle\inte\psi^2_2+\displaystyle 2\inte w\psi_2^3\Big]\div 2\pi \\[4mm]
  &=& \displaystyle\intr r(w+rw')^2-\intr rw(w+rw')^3:=A-B.
  \end{array}\end{equation}
Here we have
\[ \arraycolsep=1.5pt\begin{array}{lll}
  A=\displaystyle\intr r(w+rw')^2&=&\displaystyle\intr r^3(w')^2+\displaystyle\intr rw^2+\displaystyle\intr r^2dw^2\\[4mm]
  &=&
  \displaystyle\intr r^3(w')^2-\frac 1 2 \displaystyle\intr rw^4 ,
  \end{array} \]
where (\ref{00a:21}) is used, and
\[\arraycolsep=1.5pt\begin{array}{lll}
  B&=&  \displaystyle\intr rw(w+rw')^3\\[4mm]
  &=&\displaystyle\Big[\intr rw^4+3\displaystyle\intr r^2w^3w'\Big]+3\intr r^3w^2w'w'+\displaystyle\intr r^4w(w')^3\\[4mm]
  &=&-\displaystyle \frac 1 2 \intr rw^4 +3\displaystyle\intr r^3w^2w'w'+\displaystyle\intr r^4w(w')^3.
 \end{array}\]
Therefore, we get from (\ref{0:21}) that
\begin{equation}\label{0:22}
  \displaystyle\frac{4I}{3}\div 2\pi=\intr r^3(w')^2-\intr r^4w(w')^3-3\intr r^3w^2w'w':=C+D+E.
\end{equation}

To further simplify $I$, recall that
\begin{equation}\label{0:23}
rw''=-w'+rw-rw^3,
\end{equation}
by which we then have
\[\arraycolsep=1.5pt\begin{array}{lll}
C&=&\displaystyle\intr r^3w'dw=-\displaystyle\intr w(r^3w')'\\[4mm]
&=&-\displaystyle\intr w\big[3r^2w'+r^2(-w'+rw-rw^3)\big]\\[4mm]
&=&-\displaystyle\intr w\big[2r^2w'+r^3w-r^3w^3\big]\\[4mm]
&=&2\displaystyle\intr rw^2-\displaystyle\intr r^3w^2+\displaystyle\intr r^3w^4.
\end{array}\]
Similarly, we have
\[\arraycolsep=1.5pt\begin{array}{lll}
D=-\displaystyle\frac 1 2 \intr r^4(w')^2dw^2&=&\displaystyle\frac 1 2 \intr w^2\big[4r^3(w')^2+2r^3w'(-w'+rw-rw^3)\big]\\[4mm]
&=&  \displaystyle\intr r^3w^2(w')^2+ \displaystyle \frac 1 4 \intr r^4dw^4-\displaystyle \frac 1 6 \intr r^4dw^6\\[4mm]
&=&  \displaystyle\intr r^3w^2(w')^2- \displaystyle\intr r^3w^4+\displaystyle\frac 2 3 \intr r^3w^6.
\end{array}\]
Note from (\ref{0:23}) that
\[\arraycolsep=1.5pt\begin{array}{lll}
 -2 \displaystyle\intr r^3w^2(w')^2&=&- \displaystyle \frac 2 3   \intr r^3w'dw^3= \displaystyle\frac 2 3 \intr w^3\big[3r^2w'+r^2(-w'+rw-rw^3)\big]\\[4mm]
&=&  \displaystyle\frac 4 3 \intr w^3r^2w' + \displaystyle\frac 2 3 \intr r^3w^4- \displaystyle\frac 2 3 \intr r^3w^6\\[4mm]
&=&- \displaystyle\frac 2 3 \intr rw^4  + \displaystyle\frac 2 3 \intr r^3w^4- \displaystyle\frac 2 3 \intr r^3w^6.
\end{array}\]
We thus derive that
\[\arraycolsep=1.5pt\begin{array}{lll}
D+E&=& -2\displaystyle\intr r^3w^2(w')^2- \displaystyle\intr r^3w^4+\displaystyle\frac 2 3\intr r^3w^6\\[4mm]
&=&-\displaystyle\frac 2 3 \intr rw^4- \displaystyle\frac 1 3\intr r^3w^4 ,
\end{array}\]
by which we conclude from (\ref{00a:21}) and (\ref{0:22}) that
\begin{equation}\label{0:24}
\displaystyle\frac{4I}{3}\div 2\pi=C+D+E=\frac 1 3 \Big[2\intr rw^2-3\intr r^3w^2+2\intr r^3w^4\Big].
\end{equation}

In the following, we note that $w$ satisfies
 \begin{equation}\label{0:25}
(rw')'=rw-rw^3,\, \ r>0.
\end{equation}
Multiplying (\ref{0:25}) by $r^3 w'$ and integrating on $[0,\infty)$, we get that
\[\arraycolsep=1.5pt\begin{array}{lll}
\displaystyle\intr r^3 w'(rw')'&=&\displaystyle\intr r^3 w'[rw-rw^3]=\displaystyle\frac{1}{2}\intr r^4dw^2-\displaystyle\frac{1}{4}\intr r^4dw^4\\[4mm]
&=&-2\displaystyle\intr r^3w^2+\displaystyle\intr r^3w^4.
\end{array}\]
Note also that
\[
\intr r^3 w'(rw')'=\intr r^3 (w')^2+\frac{1}{2}\intr r^4 d(w')^2=-\intr r^3 (w')^2.
\]
By combining above two identities, it yields that
 \begin{equation}\label{0:26}\intr r^3 (w')^2 = 2\intr r^3 w^2 -\intr r^3 w^4.
 \end{equation}
On the other hand, multiplying (\ref{0:25}) by $r^2w$ and integrating on $[0,\infty)$, we obtain that
\[\arraycolsep=1.5pt\begin{array}{lll}
\displaystyle\intr r^3w^2-\displaystyle\intr r^3w^4 &=& \displaystyle\intr r^2  w(rw')'=\displaystyle\intr r^2  ww'+\displaystyle\intr r^3  wdw'\\[4mm]
&=&\displaystyle\intr r^2  ww'-\displaystyle\intr w'(3r^2w+r^3w')\\[4mm]
&=&-2\displaystyle\intr r^2  ww'-\displaystyle\intr r^3 (w')^2\\[4mm]
&=&2\displaystyle\intr rw^2-\displaystyle\intr r^3 (w')^2,
\end{array} \]
which then implies that
\begin{equation}\label{0:27}
\displaystyle\intr r^3 (w')^2=2\displaystyle\intr rw^2 -\displaystyle\intr r^3w^2+\displaystyle\intr r^3w^4.
\end{equation}
We thus conclude from (\ref{0:26}) and (\ref{0:27}) that
\[
2\intr rw^2-3\intr r^3w^2+2\intr r^3w^4=0,
\]
which therefore implies that $I=0$ in view of (\ref{0:24}), i.e., (\ref{3a:37aa}) holds.
\qed

 \vskip 0.1truein

\noindent{\em The proof of (\ref{a:claim3}).}
Following Lemmas \ref{lem2.3} and  \ref{lem3.2} again, we get that
\begin{equation}\label{A.8}
\arraycolsep=1.5pt\begin{array}{lll}
 II&=&2\displaystyle\inte  \psi_5\mathcal{L}\psi_{2}+2\displaystyle\inte \psi_1\psi_2=2\displaystyle\inte \psi_2[ \mathcal{L}\psi_5+\psi_1] \\[4mm]
  &=& -\displaystyle\inte (w+x\cdot \nabla w)(6w\psi_1\psi_2+2\psi_1) \\[4mm]
  &&-\displaystyle\frac 1 2 \inte \Big[\frac{3w^2}{a^*}+\frac{h(x+y_0)}{\lam ^{2+p}}\Big]\big(w+x\cdot \nabla w\big)^2
 \\[4mm]
  && +\displaystyle\frac{1}{2\lam ^{2+p}}\inte  (w+x\cdot \nabla w)\big[y_0\cdot \nabla h(x+y_0)\big]w\\[3mm]
  &:=&A+B.
  \end{array}\end{equation}
Since $(x+y_0)\cdot \nabla h(x+y_0)=ph(x+y_0)$ holds in $\R^2$, we derive from (\ref{1:H}) and (\ref{00:21}) that
\[
\arraycolsep=1.5pt\begin{array}{lll}
B&=&-\displaystyle\frac 1 2 \inte \Big[\frac{3w^2}{a^*}+\frac{h(x+y_0)}{\lam ^{2+p}}\Big] \big [w^2+2w(x\cdot \nabla w)+(x\cdot \nabla w)^2\big]\\[3mm]
&&+\displaystyle\frac{1}{2\lam ^{2+p}}\inte  (w+x\cdot \nabla w)\big[y_0\cdot \nabla h(x+y_0)\big]w\\[4mm]
  &=& -\displaystyle  \inte \Big[\frac{3w^2}{\inte w^4}+\frac{h(x+y_0)}{p\inte h(x+y_0)w^2}\Big] \big [w^2+2w(x\cdot \nabla w)\big]\\[4mm]
   && -\displaystyle\frac 1 2 \inte \Big[\frac{3w^2}{a^*}+\frac{h(x+y_0)}{\lam ^{2+p}}\Big] \big (x\cdot \nabla w\big)^2+\displaystyle\frac{1}{2\lam ^{2+p}}\inte w\big[y_0\cdot \nabla h(x+y_0)\big]( x\cdot \nabla w)\\[4mm]
  &=& -3-\displaystyle \frac 1 p-\displaystyle \frac{3}{2\inte w^4}\inte (x\cdot \nabla w^4)-\displaystyle \frac{1}{p\inte h(x+y_0)w^2}\inte h(x+y_0)\big(x\cdot \nabla w^2\big)\\[4mm]
   && -\displaystyle\frac 1 2 \inte \Big[\frac{3w^2}{a^*}+\frac{h(x+y_0)}{\lam ^{2+p}}\Big] \big (x\cdot \nabla w\big)^2+\displaystyle\frac{1}{2\lam ^{2+p}}\inte w\big[y_0\cdot \nabla h(x+y_0)\big]( x\cdot \nabla w)\\[4mm]
  &:=&  -3-\displaystyle \frac 1 p+3+\displaystyle \frac{2+p}{p } +C_0=\frac{p+1}{p}+C_0,
\end{array}\]
where the term $C_0$ satisfies
\[
\arraycolsep=1.5pt\begin{array}{lll}
C_0&=&-\displaystyle\frac 1 2 \inte \Big[\frac{3w^2}{a^*}+\frac{h(x+y_0)}{\lam ^{2+p}}\Big] \big (x\cdot \nabla w\big)^2+\displaystyle\frac{1}{2\lam ^{2+p}}\inte w\big[y_0\cdot \nabla h(x+y_0)\big]( x\cdot \nabla w)\\[4mm]
 &=&-\displaystyle\frac 1 {2a^*} \inte  (x\cdot \nabla w)(x\cdot \nabla w^3)-\displaystyle\frac 1 {2\lam ^{2+p} } \inte  h(x+y_0)\big(x\cdot \nabla w\big)(x\cdot \nabla w)\\[3mm]
 &&+\displaystyle\frac{1}{2\lam ^{2+p}}\inte w\big[y_0\cdot \nabla h(x+y_0)\big]( x\cdot \nabla w) \\[4mm]
  &=& \displaystyle\frac 1 {2a^*} \inte w^3\Big[ 2(x\cdot \nabla w)+x\cdot \nabla (x\cdot \nabla w)\Big]\\[4mm]
  &&+\displaystyle\frac 1 {2\lam ^{2+p} } \inte w\Big\{2h(x+y_0)(x\cdot \nabla w)+ \big[x\cdot \nabla h(x+y_0)\big]( x\cdot \nabla w)\\[3mm]
 &&\qquad\qquad +h(x+y_0)\Big[x\cdot \nabla (x\cdot \nabla w)\Big]\Big\}+\displaystyle\frac{1}{2\lam ^{2+p}}\inte w\big[y_0\cdot \nabla h(x+y_0)\big]( x\cdot \nabla w)\\[4mm]
  &=& \displaystyle\frac 1 2 \inte \Big[\frac{ w^3}{a^*}+\frac{wh(x+y_0)}{\lam ^{2+p}}\Big] \big[x\cdot \nabla \big (x\cdot \nabla w\big)\big]\\[4mm]
  &&+ \displaystyle\frac 1 {a^*} \inte w^3(x\cdot \nabla w)+\displaystyle\frac {2+p} {2\lam ^{2+p} } \inte w h(x+y_0)(x\cdot \nabla w)\\[4mm]
  &=& \displaystyle\frac 1 2 \inte \Big[\frac{ w^3}{a^*}+\frac{wh(x+y_0)}{\lam ^{2+p}}\Big] \big[x\cdot \nabla \big (x\cdot \nabla w\big)\big]\\[4mm]
  &&+ \displaystyle\frac 1 {2\inte w^4}\inte (x\cdot \nabla w^4)+\displaystyle\frac {2+p} {2p\inte h(x+y_0)w^2}\inte h(x+y_0)\big(x\cdot \nabla w^2\big)\\[4mm]
  &=& \displaystyle\frac 1 2 \inte \Big[\frac{ w^3}{a^*}+\frac{wh(x+y_0)}{\lam ^{2+p}}\Big] \big[x\cdot \nabla \big (x\cdot \nabla w\big)\big]-1-\displaystyle\frac {(2+p)^2} {2p},
\end{array}\]
in view of (\ref{00:21}).  We thus have
\begin{equation}\label{a.9}
B=\displaystyle\frac 1 2 \inte \Big[\frac{ w^3}{a^*}+\frac{wh(x+y_0)}{\lam ^{2+p}}\Big] \big[x\cdot \nabla \big (x\cdot \nabla w\big)\big]-\displaystyle\frac {p^2+4p+2} {2p}.
\end{equation}

We next calculate the term $A$ as follows. Observe that
\[
\arraycolsep=1.5pt\begin{array}{lll}
&&\displaystyle\frac 1 2 \inte \psi_1x\cdot \nabla \big (x\cdot \nabla w\big)\\[4mm]
&=&-\displaystyle\frac 1 2 \inte \big(x\cdot \nabla w\big)\big[2\psi_1+\big(x\cdot \nabla \psi_1\big)\big]\\[4mm]
  &=& -\displaystyle\inte \psi_1\big(x\cdot \nabla w\big)-\displaystyle\frac 1 2 \inte \big(x\cdot \nabla \psi_1\big)\big(x\cdot \nabla w\big)\\[4mm]
  &=& -\displaystyle\inte \psi_1\big(x\cdot \nabla w\big)+\displaystyle\frac 1 2 \inte w\big[2\big(x\cdot \nabla \psi_1\big)+
  x\cdot \nabla \big (x\cdot \nabla \psi_1\big)\big]\\[4mm]
  &=& -\displaystyle\inte \psi_1\big(x\cdot \nabla w\big)+\displaystyle\inte w\big(x\cdot \nabla \psi_1\big)+ \displaystyle\frac 1 2 \inte
 w x\cdot \nabla \big (x\cdot \nabla \psi_1\big),
\end{array}\]
which implies that
\begin{equation}\label{a.12}\arraycolsep=1.5pt\begin{array}{lll}
 &&-\displaystyle\inte \psi_1\big(x\cdot \nabla w\big)+\displaystyle\inte w\big(x\cdot \nabla \psi_1\big)\\[4mm]
 &=&\displaystyle\frac 1 2 \inte \psi_1x\cdot \nabla \big (x\cdot \nabla w\big)-\displaystyle\frac 1 2 \inte
 w x\cdot \nabla \big (x\cdot \nabla \psi_1\big).
\end{array}\end{equation}
Using (\ref{a.12}), we then derive that
\begin{equation}\label{a.10}
\arraycolsep=1.5pt\begin{array}{lll}
A&=&-\displaystyle\inte (w+x\cdot \nabla w)(6w\psi_1\psi_2+2\psi_1)  \\[4mm]
  &=& -2\displaystyle\inte w\psi_1-2\displaystyle\inte \psi_1(x\cdot \nabla w)+3\displaystyle\inte w\psi_1(w+x\cdot \nabla w)^2\\[4mm]
  &=& -2\displaystyle\inte w\psi_1- \displaystyle\inte \psi_1(x\cdot \nabla w)\\[4mm]
  &&+\displaystyle\inte w\big[2\psi_1+x\cdot \nabla \psi_1\big]+3\displaystyle\inte w\psi_1(w+x\cdot \nabla w)^2\\[4mm]
  &=&  - \displaystyle\inte \psi_1(x\cdot \nabla w)+\displaystyle\inte w(x\cdot \nabla \psi_1)+D\\[4mm]
  &=& \displaystyle\frac 1 2 \inte \psi_1x\cdot \nabla \big (x\cdot \nabla w\big)-\displaystyle\frac 1 2 \inte
 w x\cdot \nabla \big (x\cdot \nabla \psi_1\big)+D,
  \end{array}\end{equation}
where the term $D$ satisfies
\[\arraycolsep=1.5pt\begin{array}{lll}
D&=&3\displaystyle\inte w\psi_1\Big[w^2+2w(x\cdot \nabla w)+(x\cdot \nabla w)^2\Big]\\[4mm]
  &=&3\displaystyle\inte w^3\psi_1+6\displaystyle\inte w^2\psi_1(x\cdot \nabla w)+\displaystyle\frac 3 2 \inte \psi_1(x\cdot \nabla w)(x\cdot \nabla w^2)\\[4mm]
  &=&3\displaystyle\inte w^3\psi_1+6\displaystyle\inte w^2\psi_1(x\cdot \nabla w)\\[4mm]
 && -\displaystyle\frac 3 2 \inte w^2\Big\{2\psi_1(x\cdot \nabla w)+(x\cdot \nabla w)(x\cdot \nabla \psi_1)+
 \psi_1\Big[x\cdot \nabla (x\cdot \nabla w)\Big]\Big\}.
\end{array}\]
Since
\[\arraycolsep=1.5pt\begin{array}{lll}
&&-\displaystyle\frac 3 2 \inte w^2(x\cdot \nabla w)(x\cdot \nabla \psi_1)\\[4mm]
&=&-\displaystyle\frac 1 2 \inte (x\cdot \nabla \psi_1)(x\cdot \nabla w^3)\\[4mm]
&=&\displaystyle\frac 1 2 \inte w^3\Big[x\cdot \nabla (x\cdot \nabla \psi_1)+2(x\cdot \nabla \psi_1)\Big]\\[4mm]
  &=&\displaystyle\frac 1 2 \inte w^3x\cdot \nabla (x\cdot \nabla \psi_1)+\displaystyle  \inte w^3(x\cdot \nabla \psi_1)\\[4mm]
  &=&\displaystyle\frac 1 2 \inte w^3x\cdot \nabla (x\cdot \nabla \psi_1)-\displaystyle  \inte \psi_1\Big[2w^3+3w^2(x\cdot \nabla w)\Big]\\[4mm]
  &=&\displaystyle\frac 1 2 \inte w^3x\cdot \nabla (x\cdot \nabla \psi_1)-2\displaystyle  \inte w^3\psi_1-3 \displaystyle  \inte w^2\psi_1(x\cdot \nabla w),
\end{array}\]
the term $D$ can be further simplified as
\begin{equation}\label{a.11}
D= \displaystyle\inte w^3\psi_1-\displaystyle\frac 3 2 \inte w^2 \psi_1\big[x\cdot \nabla (x\cdot \nabla w)\big]+\displaystyle\frac 1 2 \inte w^3x\cdot \nabla (x\cdot \nabla \psi_1).
\end{equation}
Applying  (\ref{a.11}), we then obtain from  (\ref{a.10}) that
\begin{equation}\label{a.14}\arraycolsep=1.5pt\begin{array}{lll}
A&=& \displaystyle\inte w^3\psi_1+\displaystyle\frac 1 2 \inte (1-3w^2) \psi_1\big[x\cdot \nabla (x\cdot \nabla w)\big]\\[4mm]
&&-\displaystyle\frac 1 2 \inte \Delta w \big[x\cdot \nabla (x\cdot \nabla \psi_1)\big],
\end{array}\end{equation}
since $w$ solves the equation $w^3-w=-\Delta w$ in $\R^2$.

Combining (\ref{a.9}) and (\ref{a.14}) now yields that
\begin{equation}\label{a.15}\arraycolsep=1.5pt\begin{array}{lll}
II=A+B&=& \displaystyle\inte w^3\psi_1-\displaystyle\frac {p^2+4p+2} {2p}\\[4mm]
&&+\displaystyle\frac 1 2 \inte \big[x\cdot \nabla (x\cdot \nabla w)\big]\Delta  \psi_1-\displaystyle\frac 1 2 \inte \big[x\cdot \nabla (x\cdot \nabla \psi_1)\big]\Delta w.
\end{array}\end{equation}
We claim that
\begin{equation}\label{a.16}
 \inte w^3\psi_1=\frac{p+1}{p}.
\end{equation}
Actually, multiplying (\ref{00:21}) by $w$ and integrating on $\R^2$ gives that
\[
\inte \nabla \psi_1\nabla w-3 \inte w^3\psi_1=-\inte \Big[\frac{2w^4}{\inte w^4}+\frac{2h(x+y_0)w^2}{p\inte h(x+y_0)w^2}\Big]=-\frac{2(p+1)}{p},
\]
due to the fact that $\inte w\psi_1=0$ by (\ref{3a:37a}). On the other hand, multiplying (\ref{Kwong}) by $\psi_1$ and integrating on $\R^2$ gives that
\[
\inte \nabla \psi_1\nabla w=-\inte w\psi_1+ \inte w^3\psi_1= \inte w^3\psi_1.
\]
The claim (\ref{a.16}) then follows directly  from above two identities. We next claim that
\begin{equation}\label{a.17}
\inte \big[x\cdot \nabla (x\cdot \nabla w)\big]\Delta  \psi_1=
\inte \big[x\cdot \nabla (x\cdot \nabla \psi_1)\big]\Delta w.
\end{equation}
To prove (\ref{a.17}), rewrite $\psi_1$ as $\psi_1(x)=\psi_1(r,\theta )$, where $(r,\theta )$ is the polar coordinate in $\R^2$, such that
\begin{equation}\label{a.18}
\Delta \psi_1 = \big(\psi_1\big)_{rr} +\frac{1}{r}\big(\psi_1\big)_r+ \frac{1}{r^2}\big(\psi_1\big)_ {\theta\theta},\quad    \nabla \psi_1=\frac{x}{r}\big(\psi_1\big)_r+\frac{x^\bot}{r^2}\big(\psi_1\big)_\theta,
\end{equation}
where $x^\bot =(-x_2,x_1)$ for $x=(x_1,x_2)\in \R^2$. We then derive from (\ref{lem2.1:2}) that
\[\arraycolsep=1.5pt\begin{array}{lll}
\displaystyle\inte \big[x\cdot \nabla (x\cdot \nabla w)\big]\Delta  \psi_1&=&\displaystyle\int ^{2\pi}_0\intr r(rw')'\Big\{\big[r\big(\psi_1\big)_r\big]_r+\frac{\big(\psi_1\big)_{\theta\theta}}{r}\Big\}drd\theta\\[4mm]
&=&\displaystyle\int ^{2\pi}_0\intr r(rw')'\big[r\big(\psi_1\big)_r\big]_rdrd\theta+\displaystyle\int ^{2\pi}_0\intr (rw')'\big(\psi_1\big)_{\theta\theta}drd\theta\\[4mm]
&=&\displaystyle\int ^{2\pi}_0\intr r(rw')'\big[r\big(\psi_1\big)_r\big]_rdrd\theta,
\end{array}\]
and
\[\inte \big[x\cdot \nabla (x\cdot \nabla \psi_1)\big]\Delta w=\displaystyle\int ^{2\pi}_0\intr r\big[r\big(\psi_1\big)_r\big]_r(rw')'drd\theta,
\]
which thus imply that (\ref{a.17}) holds. Applying (\ref{a.16}) and (\ref{a.17}), we therefore conclude from (\ref{a.15}) that
\[
II=\frac{p+1}{p}-\frac {p^2+4p+2} {2p}=-\frac{2+p}{2},
\]
which gives (\ref{a:claim3}), and the proof is complete.
\qed
 \vskip 0.1truein

 \vskip 0.16truein
\noindent {\bf Acknowledgements:} The authors are  very grateful to the referees for many useful suggestions which lead to some improvements of the present paper.  The first author thanks Prof. Robert Seiringer for fruitful discussions on the present work. Part of the present work was finished when the first author was visiting  Taida Institute of Mathematical Sciences (TIMS) in October 2013 and Pacific Institute for Mathematical Sciences (PIMS) from March to April in
2016. He would like to thank both institutes for their warm hospitality.



\end{document}